\documentclass[a4paper]{article}
\usepackage{amsmath,amssymb,amsthm,amsfonts,epsfig,fullpage}
\usepackage{color}
\usepackage[applemac]{inputenc}

\def\Rset{\mathbb{R}}
\def\Nset{\mathbb{N}}

\def\Mset{\mathbb{M}}

\def\neig{{\mathrm {Neighbors}}}

\def\diag{{\mathrm {\tt diag}}}
\def\diam{{\mathrm {\tt diam}}}

\def\bfo{{\mathbf 1}}

\def\cA{{\cal A}}

\def\cI{{\cal I}}
\def\cL{{\cal L}}
\def\cM{{\cal M}}
\def\cN{{\cal N}}

\def\cT{{\cal T}}

\def\trho{\tilde{\rho}}

\def\neig{{\mathrm {Neighbors}}}   
   
\def\card{{\mathrm {card}}}

\def\cN{{\cal G}}
\def\cI{{\cal I}}
\def\cN{{\cal N}}

\newcommand{\comment}[1]{}
\renewcommand{\t}{^{\mbox{\tiny\sf T}}}

\definecolor{extra}{rgb}{0.0, 0.75, 0.55}

\newtheorem{theo}{Theorem}
\newtheorem*{theo*}{Theorem}
\newtheorem{defi}{Definition}

\newtheorem{lemm}[theo]{Lemma}
\newtheorem{coro}[theo]{Corollary}

\newtheorem*{prop*}{Proposition}
\newtheorem{prop}[theo]{Proposition}

\theoremstyle{definition}
\newtheorem{rema}{Remark}

\newtheorem{exam}{Example}

\begin{document}
\title{Convergence speed of unsteady distributed consensus:
decay estimate along the settling spanning-trees}
\author{David Angeli\thanks{Dipartimento di Sistemi e Informatica, 
University of Florence, Via di S.\ Marta 3, 50139 Firenze, Italy. 
Email: 
\tt angeli@dsi.unifi.it} \and Pierre-Alexandre Bliman\thanks{INRIA, 
Rocquencourt BP105, 78153 Le Chesnay cedex, France. Email: \tt 
pierre-alexandre.bliman@inria.fr}}
\date{\today}
\maketitle

\begin{abstract}
Results for estimating the convergence rate of non-stationary distributed consensus
algorithms are provided,
on the basis of qualitative (mainly topological) as well as basic quantitative information (lower-bounds on the 
matrix entries). The results appear to be tight in a number of instances and are illustrated through simple as well
as more sophisticated examples. The main idea is to follow propagation of information along certain spanning trees which arise in the communication graph.
\end{abstract}

\section{Introduction}
\label{se1}

Historically appeared in the areas of communication networks, control theory and parallel computation, the analytical study of ways for reaching consensus in a population of agents is a problem of broad
interest in many fields of science and technology.
Questions of this nature arise in peer-to-peer 
and sensors networks \cite{boyd,akar}, in the manoeuvring of groups of vehicles \cite{AND,Naomi,YAN}, in the study of TCP protocols \cite{shorten}, in the theory of coupled oscillators \cite{GAD,MAN99a,JadMot2004,Sepulchre}, in neural networks \cite{Hirsch}, but also in apparently distant fields such as in the study and modeling of opinion in social science \cite{lorenz} or of animal flocking \cite{CZI}.

Generally speaking, the aim of such studies is to design or analyse decentralised algorithms through which ``agents" (which in the previous examples can be cars or unmanned aerial vehicles, nodes in communication network, sensors, particles, cells, fish\dots) can update their internal states
in order to agree on a common value for such variable.
In general, the latter shall not be a priori fixed, but will be determined as a result of the interactions and of their history\footnote{Distributed averaging corresponds to the special case where the limit is guaranteed to be the average of the individual values.}.
These interactions can be modeled either as unidirectional or bidirectional, corresponding to different extreme
situations in which one agent is able to influence somebody else without being affected by the internal
state of the receiver (as in hierarchical communication flows) or,  the opposite situation,
in which influence between agents is always symmetrical.

Of particular interest is the question of estimating how quickly consensus is
reached on the basis of few qualitative (mainly topological) information as well as basic quantitative  
(mainly the strength of reciprocal influences) information on the network.    

Originally, the problem of quantifying the convergence rate towards consensus was considered mainly in the context of
stationary networks. For Markov chains, i.e., this amounts to quantify the speed at which 
steady-state probability distribution is achieved, and is therefore directly related to finding an
a priori estimate to the second largest eigenvalue of a stochastic matrix.
Classical works on this subject are due to Cheeger and Diaconis, \cite{cheeger,diaconis}, see also
\cite{friedland} for improved bounds.

Among the classical contributions which instead deal with time-varying interactions we refer to the work of Cohn, \cite{cohn}, where asymptotic convergence is proved, but neglecting the issue of relating topology and
guaranteed convergence rates.
J.N.\ Tsitsiklis {\em et al.} also provided important qualitative contributions to this subject \cite{TSI84,TSI86,BER}, as well as L.\ Moreau \cite{MOR04a}.
See also \cite{ANG} for further nonlinear results.
In particular, the role of connectivity of the communication graph in the convergence of consensus has been recognised and finely analysed.


As noticed in different manners by the preceding authors, arguments based on graph theory are more powerful and seem to catch in a more natural way the essence of the problem, rather than computations based on linear algebra techniques (although the study of stochastic matrices offer nowadays, undoubtedly, quite strong results).
We are in perfect harmony with the opinion that vision in terms of graph is central to understand the agreement issues.
However, it appears that some dynamical aspects which have been so far disregarded can be exploited
to really gain a tighter understanding of how rapidly consensus can be reached.
Our attempt here is to provide a consequent step toward integration of the temporal aspects of information transit.
We are thus led to further elaborate and exploit tools for description of the connectivity emergence in the communication graphs.\\

Our purpose in this paper is to provide several criteria to estimate quantitatively the contraction rate of a set of agents towards consensus, in a discrete time framework.
Using the language of dynamical system, the problem is here of estimating the second largest Lyapunov exponent of an infinite product of matrices (see also \cite{BLO} for links with some joint spectral radius).
To the best of our knowledge, previous results are centrally based on the existence of a lower bound of the nonzero entries
associated to such matrices, and most of them on the existence of self-loops, see \cite{CAO05} and the surveys in \cite{BLO,OLS} (see however the contributions in \cite{TSI84,TSI86} where the assumption on self-loops is relaxed).
 Recently, A.\ Nedich {\em et al.} \cite{NED} proposed improved bounds under similar assumptions.
On the contrary, we attempt here to follow more closely the spread of the information over the agent population, along the one or more spanning-trees.

Ensuring a lower bound to the matrix entries of the agents already attained by the information flow along the spanning-tree, rather than all the nonzero contributions as classically, permits to obtain tighter estimates with weaker assumptions. 
The setting used here applies indifferently to leader-follower or to leaderless networks.



More precisely, the main idea is to examine the birth and rise of spanning-trees in the network.
Distinguishing between different sub-populations, of agents already touched by spanning-tree and of agents not yet attained, and using lower bounds on the influence of the former ones on the latter ones, one is able to establish rather precise convergence estimates.
Due to the nature of the assumptions, the latter possess some innate robustness with respect to parametric uncertainties.

The paper is organized as follows.
The problem is formulated in Section \ref{se3}, and some pertinent concepts are therein introduced.
Specifically, several appropriate connectivity notions are defined, among which {\em sequential connectivity}, which turns out to be central to our developments.
The remaining Sections are devoted to the statement and demonstration of the main results.
 Section \ref{se4} deals with the problem of contraction estimates when information follows a single spanning-tree (or at least one such spanning-tree with a certain guaranteed strength exists in the underlying graph).
 This section contains remarks on self-loops and delays, see Subsection \ref{se33}.
The most original results are in Section \ref{se55}, where it is shown by means of a fairly general technique, how multiple spanning-trees can be used to derive tight estimates of the contraction rate.
Last, conclusions are reported in Section \ref{se77}.

For better readability, various examples are reported in the text, to illustrate the application of the results and to demonstrate the powerfulness of the method.
Also, some involved results and proofs  have been put in Appendix.
In particular the main technical tool for carrying out estimates over a finite horizon of the contraction rate of a linear stochastic system is given there.\\

The authors wish to thank the reviewers and the Associate Editor for their penetrating comments, including the detection of a mistake in the previous version.





\paragraph{Notation}

In the sequel, $\Nset$ stands for the set of natural integers (including zero), $\lfloor x \rfloor$ designates the integer value of a real number $x$.
For any set $\cN$, we denote $|\cN|$ or $\card\ \cN$ its cardinality.
Generally speaking (Latin or Greek) upper case letters indicate matrices, and lower case letters are used to signal scalar numbers and vectors.
Graphs and sets are distinguished by calligraphic letters.

 In the sequel, the (possibly infinite-dimentional) vectors $\bfo$ and $\bfo_i$ denote respectively a column of $1$ and the vector with null components, except 1 in the $i$-th position.

We call {\em integer interval} any set obtained as the intersection of a usual interval with the set $\Nset$.
When the context is clear, in particular when talking about time values, the integer intervals are denoted as the classical ones: for example, $[0,T]\doteq\{t\in\Nset\ :\ 0\leq t\leq T\}$.

For $p,q$ positive integers, we denote as usual $I_p$ and $0_{p\times q}$, the identity and zero matrices.
The transposition of matrices is denoted $\t$.
By definition, (row) stochastic (resp.\ sub-stochastic) matrices are square matrices with nonnegative components, whose row sums are equal (resp.\ at most equal) to 1.
Their spectrum is ordered by nonincreasing modulus magnitude: for $M$ stochastic in $\Rset^{n\times n}$, $1=\lambda_1(M)\geq |\lambda_2(M)|\geq\dots |\lambda_n(M)|$.

Last, we introduce the matrix sets $\Mset^{p, q}$.
By definition,
\begin{equation}
\label{Mpq}
\Mset^{p, q}
\doteq
\left\{
M\in\Rset^{p\times q}\ :\ M\geq 0 \text{ and } \forall i=1,\dots, p,\ M_{i,1}+\dots+M_{i,q} \geq 1
\right\}\ .
\end{equation}
In \eqref{Mpq} and everywhere in the paper, matrix ordering is meant componentwise: $M\geq 0$ stands for $M_{i,j}\geq 0$ for all $i,j$. 
 
\section{ Sequential connectivity and other graph related notions}
\label{se3}

We consider the problem of convergence of the consensus algorithm described by the following system:
\begin{equation}
\label{equb}
x_k(t+1) = \sum_{l\in\cN} \gamma_{k,l}(t)x_l(t),\ k\in\cN\ ,
\end{equation}
toward a {\em common value}; that is, the {\em global asymptotic stability of the diagonal set} $\{x\ :\ \forall k,l\in\cN, x_k=x_l\}$.
As usual, \eqref{equb} may be written in matrix form, as
\[
x(t+1)=\Gamma(t)x(t),\quad
x(t)\doteq\begin{pmatrix} x_1(t)\\ \vdots \end{pmatrix},\
\Gamma(t)\doteq (\gamma_{k,l}(t))_{(k,l)\in\cN\times\cN}\ .
\]
We consider scalar systems, although extension to multidimensional systems is possible.
The set $\cN$ is finite or countable, and the functions $x_k$ map $\Nset$ to $\Rset$.
We assume in all the sequel that
\begin{equation}
\label{alpha}
\forall k,l\in\cN,\forall t\in\Nset,\ \gamma_{k,l}(t)\geq 0,
\text{ and } \forall k\in\cN,\ \sum_{l\in\cN}\gamma_{k,l}(t)=1\ .
\end{equation}
In other words, the matrices $(\gamma_{k,l}(t))_{(k,l)\in\cN\times\cN}$ are stochastic.

Our goal in the remaining of the paper is to quantify the convergence speed of the set $\{x_k(t)\ :\ k\in\cN\}$ when $t\to +\infty$ toward a consensus value.
We first introduce vocabulary adequate to measure the latter.

\begin{defi}[Agent set diameter]
\label{de1}
The quantity
\[
\Delta(x(t)) \doteq \sup_{k\in\cN} x_k(t) - \inf_{k\in\cN} x_k(t)
\]
is called the {\em diameter of the agent set at time $t$.}
\end{defi}
In what follows, $\Delta(x(t))$ plays the role of a Lyapunov function to study convergence to an agreement.
Although the latter depends upon the state, we frequently abbreviate the notation in $\Delta(t)$ if no misinterpretation is possible.
\begin{defi}[Contraction rate]
\label{de15}
We call {\em contraction rate\/} of system \eqref{equb} the number $\rho\in [0,+\infty]$ defined as:
\[
\rho \doteq \sup_{x(0)}\
\limsup_{t\to +\infty} \left ( \frac{\Delta(t)}{\Delta(0)}  \right )^{\frac{1}{t}}\ .
\] 
\end{defi}
The number $\rho$ is indeed the second largest Lyapunov exponent of the dynamical system \eqref{equb}.

Some notions and definitions necessary to describe pertinent aspects of the communication between the agents are now introduced, based on some elementary tools of algebraic graph theory.

\begin{defi}[Communication graph]
\label{de2}
We call {\em communication graph} (of system \eqref{equb}) {\em at time $t$}
the directed graph defined by the ordered pairs $(k,l)\in\cN\times\cN$ such that $\gamma_{k,l}(t)>0$.
\end{defi}
In the present context, we use indifferently the terms ``node" or ``agent".

\begin{defi}[Neighbors]
\label{de25}
Given a graph $\cA$ and a nonempty subset 
$\cL\subseteq\cN$, the set $\neig(\cL,\cA)$ of {\em neighbors of $\cL$} is the set of those agents $k\in\cN\setminus\cL$ for which there exists at least one element $l\in\cL$ such that $(k,l)\in\cA$.
When $\cL$ is a singleton $\{l\}$, the notation $\neig(l,\cA)$ is used 
instead of $\neig(\{l\},\cA)$.

\end{defi}
A key property, namely {\em weak connectivity}, has been shown to influence crucially the evolution of finite systems of agents linked by time-varying communication graphs (see \cite{MOR03,MOR04a}, but also \cite{CAO06}, where the weakly connected sequences are called ``repeatedly jointly rooted").

\begin{defi}[Connectivity and weak connectivity]
\label{de3}
A node $k\in\cN$ is said to be {\em connected\/} to a node  $l\in\cN$ on a directed graph $\cA$ defined on $\cN$, if there exists a {\em path\/} joining $k$ to $l$ in $\cA$ and respecting the orientation of the arcs.
Given a sequence of directed graphs $\cA(t)$, $t\in\Nset$, the node $k\in\cN$ is said {\em connected\/} to the node $l\in\cN$ {\em on an integer interval $I\subseteq\Nset$\/} if $k$ is connected to $l$ for the graph $\bigcup_{t\in I}\cA(t)$.

A graph $\cA$ is called {\em weakly connected\/} {\rm\cite{MOR03}} if there is a node $k\in\cN$ connected to all other nodes $l \in\cN$.
A sequence of graphs $\cA(t)$, $t\in\Nset$,  is called {\em weakly connected across
an integer interval $I\subseteq\Nset$} if the graph $\bigcup_{t\in I}\cA(t)$ is weakly connected (that is, if there is a node connected across $I$ to all other nodes).
A subgraph connecting an agent to all the other ones is called a {\em spanning-tree}.
\end{defi}

The fundamental result found by Moreau states that uniform global asymptotic stability of the set of common equilibria is {\em equivalent\/} to the existence of an integer $T>0$ such that the sequence of graphs is weakly connected on any interval of length $T$ \cite{MOR03,MOR04a}.
Exponential estimates may be obtained too, see the survey part of \cite{BLO,OLS}, and \cite{CAO05,CAO06}.
As a matter of fact, there is no specific difficulty to check the validity of both these results, with the weaker assumption that the graph sequence is weakly connected on every integer intervals  $[t_p, t_{p+1}]$, $p\in\Nset$, where the $t_p$ define a strictly increasing sequence such that $\limsup_{p \to +\infty} t_{p+1}-t_p \leq T$.\\

In order to obtain more precise estimates of the decay rate toward consensus value, it is reasonable to introduce some minimal time taken by the information to cover the graph --- while the preceding connectivity notions were not concerned with the ordering of the arcs constituting the tree.
 We thus introduce in the sequel some notions useful to quantify the minimal time for information spread.
The latter play a central part in the contraction rate estimate to be stated later.

\begin{defi}[Sequential connectivity of finite graph sequences]
\label{de4bis}
A finite sequence of $T$ graphs with common nodes $\cA_1, \cA_2$, \dots, $\cA_T$ ($T\in\Nset$) is said to be {\em sequentially connected} if
there exist a node $k \in \cN$ and iterations given by:
\begin{eqnarray} 
\label{spantree}
\cN_0 &=& \{ k \} \nonumber \\
\cN_{t} &\subseteq& \cN_{t-1} \cup \neig ( \cN_{t-1}, \cA_{t} ) \quad t = 1, \ldots, T \nonumber
\end{eqnarray}
which satisfy $\cN_{T}= \cN$. \\
\end{defi}
When we want to emphasize the ``root'' node, 
we denote $\cN_t$ by $\cN_t ( k ) $, meaning that the iteration departs from node $k$. 

The sets introduced in Definition \ref{de4bis} are crucial to understand the principle of the method developed in the present paper.
For each $t=1,\dots, T$, the set $\cN_t$ contains agents already in $\cN_{t-1}$ and agents having a neighbour in $\cN_{t-1}$ at time $t$: they are all agents which have been attained at most at time $t$ by the settling of the spanning-tree rooted in $k$.

We now introduce a derived notion for infinite sequences of graphs.

\begin{defi}[$T$-sequential connectivity]
\label{de5}
An infinite sequence of graphs $\cA (t)$, $t \in \Nset$, is said {\em $T$-sequentially connected} if there 
exists a strictly increasing integer sequence $t_p$, $p\in\Nset$, fulfilling
\begin{equation}
\label{defisup}
\limsup_{p \rightarrow + \infty} \, t_{p+1} - t_p \leq T
\end{equation}
and such that  each graph sub-sequence
\[
\cA ( t_p ),  \ldots, \cA (t_{p+1}-1)
\]
is sequentially connected.
\end{defi}
Remark that the property is by definition monotone with respect to $T$, viz.
\[ \textrm{$T$-sequential connectivity} \, \Rightarrow \, \textrm{$(T+1)$-sequential connectivity}. \]
Moreover, $T$-sequential connectivity is invariant with respect to finite time shifts, namely 
$\cA(t)$ is $T$-sequentially connected iff for all $q \in \Nset$, $\cA(t+ q)$ is again
 $T$-sequentially connected.
 Similarly, $T$-sequential connectivity is invariant with respect to deletions and/or substitutions of
 finitely many graphs in a sequence, thus confirming that the property is truly an asymptotic
 definition.

Notice the proximity of the definitions of sequential connectivity proposed here, with the notion of weak connectivity;
the central difference being that the former one takes into account explicitly the time scheduling of the information transit. 
The following result links the different connectivity properties defined above and provides mutual bounds between the different connectivity time constants.
\begin{prop}
\label{pr1}
Any $T$-sequentially connected sequence of graphs is weakly connected on the integer intervals $[t_p, t_{p+1}]$, $p\in\Nset$.
Reciprocally, given an increasing sequence $t_p$ fulfilling \eqref{defisup}, any sequence of graphs defined on a set of $n$ agents that is weakly connected on the intervals $[t_p, t_{p+1}]$, $p\in\Nset$, is $(n-1)^2T$-sequentially connected.
\end{prop}

\begin{proof}
The first statement is straightforward.
We  show next the converse part. For each $p$ in $\Nset$, let $h_p$ denote any of the agents connected to all the other ones over the union of graphs $\bigcup_{t=t_p}^{t_{p+1} -1} \cA (t)$ (such an $h_p$ 
always exists because of weak connectivity).
Consider the sequence $h_p, h_{p+1}, \ldots h_{p+(n-1)^2-1}$ (of length $(n-1)^2$).   
Since $(n-1)^2=(n-1)+(n-2)(n-1)$, it becomes obvious that at least one agent appears $(n-1)$
times or more along the above sequence. Let us denote this agent by $k$ and $[t_{p_l}, t_{p_l+1}]$, $l=1,\dots, n-1$
the corresponding time intervals.
Remark that for large enough $i$, $[t_p, t_{p+(n-1)^2}]\subseteq [t_p, t_p+(n-1)^2T]$, due to \eqref{defisup}.

Define $\cN_{t_{p_l}}$, $l=0, \dots, n-1$, by $\cN_{t_{p_0}}\doteq\{k\}$ and
\begin{equation}
\label{cNp}
\cN_{t_{p_{l+1}}} \doteq \cN_{t_{p_l}}\cup\left\{
l'\in\cN\ :\ \exists t=1,\ \dots,\ t_{p_{l+1}}-t_{p_l},\ \neig(l',\cA(t_{p_l}+t))\cap\cN_{t_{p_l}}\neq\emptyset
\right\}\ .
\end{equation}
The set $\cN_{t_{p_0}}$ being a singleton, it suffices to show that the sequence $\cN_{t_{p_l}}$ is increasing, in order to deduce that $|\cN_{t_{p_{n-1}}}|=n=|\cN|$, and so $\cN_{t_{p_{n-1}}}=\cN$.
The latter property will then imply the existence of a sequential spanning-tree in $[t_{p_0}, t_{p_{n-1}}]$ (included in $[t_p, t_p+(n-1)^2T]$ for large enough $i$), and consequently the claimed proposition that the sequence of graphs is $(n-1)^2T$-sequentially connected.

Let $l\in\{0, 1,\dots , n-1\}$, let us thus establish that $\cN_{t_{p_l}}\subsetneq\cN_{t_{p_{l+1}}}$, provided that $\cN_{t_{p_l}}\neq\cN$.
By assumption, agent $k$ is connected ({\em non}-sequentially) to all other agents at least on the $n-1$ different time integer intervals $[t_{p_l}, t_{p_{l+1}}]$.
Thus, there exists $t\in \{1,\dots, t_{p_{l+1}}-t_{p_l}\}$ such that the graph $\cA(t_{p_l}+t)$ contains a link originating in $\cN_{t_{p_l}}$ and terminating outside.
This implies that, for this value of $t$, $\neig(l',\cA(t_{p_l}+t))\cap\cN_{t_{p_l}}\neq\emptyset$.
In view of \eqref{cNp}, one thus deduces that $\cN_{t_{p_l}}\subsetneq\cN_{t_{p_{l+1}}}$.
This achieves the proof of Proposition \ref{pr1}.
\end{proof}


\section{Propagation of a unique spanning-tree}
\label{se4}


\subsection{Estimating the contraction: a key lemma}
\label{se31}

A first result is now given, describing the elementary mechanism which permits to quantify contraction along a {\em unique} spanning-tree.

\begin{lemm}
\label{le0}
Let the finite sequence of communication graphs $\cA(0),\dots, \cA(T-1)$ of system \eqref{equb} be sequentially connected, and let $\cN_0,\dots,\cN_T$ be the sets corresponding to the spanning-tree (see Definition {\rm\ref{de4bis}}).
Assume that, for any $t=0,\dots, T-1$ and any $k\in\cN$,
\begin{equation}
\label{exemm}
k \in \cN_{t+1} 
\Rightarrow
\sum_{l\in\cN_{t} } 
 \gamma_{k,l}(t) \geq\alpha(t)\ ,
\end{equation}
for a given map $\alpha: [0,T-1]\to [0,1]$.
Then
\begin{equation}
\label{recaa}
\Delta(T) 
\leq \left(
1-\prod_{t=0}^{T-1}\alpha(t)
\right)\Delta(0)\ .
\end{equation}
\end{lemm}

Besides sequential connectivity, it is thus assumed in Lemma \ref{le0} that, when an agent is in $\cN_{t+1}$ (and thus is attained by the spanning-tree at most at time $t+1$), then at time $t$ the total weight in the right-hand side of \eqref{equb} of its neighbours from $\cN_{t}$ (which have been previously attained by the spanning-tree), including possibly itself, is at least $\alpha(t)$, until completion of the tree.
This is thus an hypothesis on the relative value of the two ``feeding weights", internal and external to the spanning-tree.

\begin{proof}[Proof of Lemma {\rm\ref{le0}}]
Lemma \ref{le0} is a particular case of a more complex result, Lemma \ref{le1}, which will be used and demonstrated further.
For this reason, we limit the present proof to the essential arguments.
For any $k \in \cN_{t+1}$  we have
\[
x_k (t+1) = \sum_{l=1,\dots, n} \gamma_{k,l}(t)x_l(t)
= \sum_{l\in\cN_{t}} 
 \gamma_{k,l}(t)x_l(t)
+ \sum_{l\in\cN \setminus\cN_{t}} 
\gamma_{k,l}(t)x_l(t)\ ,
\]
where by assumption $\sum_{l\in\cN_{t}} 
\gamma_{k,l}(t) \geq \alpha(t)$ and $\sum_{l=1,\dots, n} \gamma_{k,l}(t)=1$.
From this, it may be shown that
\[
\max_{k\in\cN_{t+1}} x_k(t+1)\leq \alpha (t) \, \max_{k\in\cN_t} x_k(t)+(1-\alpha(t) ) \, \max_{k\in\cN} x_k(t)\ ,
\]
and  opposite inequality for the corresponding $\min$ expressions.
Denoting
\[
\Delta_{1}(t)\doteq\max_{k\in\cN_{t}}x_k(t)-\min_{k\in\cN_{t}}x_k(t)\ ,
\]
it turns out that $\Delta_{1}(0)=0$ ($\cN_{0}$ is a singleton, the root of the spanning-tree), while $\Delta_{1}(T)=\Delta(T)$ (when the spanning-tree has run entirely the graph, at time $T$).
Thus
\[
\Delta(T)
= \Delta_{1}(T)
\leq \left( 1-\alpha(0)\alpha(1)\dots\alpha(T-1)
\right) \, \Delta(0)
\]
as claimed in the statement.
\end{proof}

\begin{rema}
\label{re2}
Under the hypotheses of Lemma \ref{le0}, one may show easily that $x(T)$, considered as a function of $x(0)$, verifies:
\begin{equation}
\label{clev1}
\forall l\in\cN,\ \frac{\partial x_l(T)}{\partial x_k(0)} \geq \prod_{t=0}^{T-1} \alpha(t)\ ,
\end{equation}
where by $k$ is denoted the index of the root of the spanning-tree.
Indeed, one has more generally, for any $t=0,\dots, T-1$:
\[
\forall l\in\cN_{t+1},\ \frac{\partial x_l(t+1)}{\partial x_k(0)} \geq \alpha(0) \alpha(1)\dots\alpha(t)\ .
\]
From \eqref{clev1}, one deduces that, for any $l\in\cN$,
\[
x_l(T)= \prod_{t=0}^{T-1}\alpha(t)\ x_k(0)+\sum_{l'\in\cN}\zeta_{l,l'}(0,T) x_{l'}(0)\ ,
\]
where the integers $\zeta_{l,l'}(0,T)$ are nonnegative and such that, for any $l\in\cN$, the weight of the remaining influences verifies:
\[
\sum_{l'\in\cN}\zeta_{l,l'}(0,T)= 1-\prod_{t=0}^{T-1}\alpha(t)\ .
\]
It is then immediate to establish that
\begin{multline*}
\left(
1- \prod_{t=0}^{T-1}\alpha (t)
\right)
\min_{l\in\cN}x_l(0)
\leq
\min_{l\in\cN}x_l(T)- \prod_{t=0}^{T-1}\alpha (t)\ x_k(0)\\
\leq
\max_{l\in\cN}x_l(T)- \prod_{t=0}^{T-1}\alpha (t)\ x_k(0)
\leq
\left(
1- \prod_{t=0}^{T-1}\alpha  (t)
\right)
\max_{l\in\cN}x_l(0)\ ,
\end{multline*}
which furnishes an alternative proof of inequality \eqref{recaa}.
One can thus interpret Assumption \eqref{exemm} as ensuring a minimal  guaranteed influence of the value of $x_k(0)$ (at the root of the spanning-tree) on every value $x_l(T)$, $l \in\cN$.
\qed
\end{rema}

\subsection{Results on contraction rate estimate}
\label{se32}

The main result of Section \ref{se4} is now presented.
Direct consequence of Lemma \ref{le0}, it provides an estimate of the contraction rate.

\begin{theo}
\label{th1}
Let the sequence of  communication graphs of system \eqref{equb} be $T$-sequentially connected.
Accordingly, denote $t_p$ the corresponding increasing  integer sequence of spanning-tree completion (see Definition {\rm\ref{de5}}); $\cN_{p,t-t_p}$, $t=t_p, t_p+1,\dots, t_{p+1}$, the sets corresponding to the spanning-tree connecting sequentially the graph sub-sequences $\cA(t_p), \dots, \cA(t_{p+1}-1)$, $p\in\Nset$ (see Definition {\rm\ref{de4bis}});
and $\cT$ the corresponding set
\begin{equation}
\label{cT}
\cT\doteq\{ (p,t)\ :\ p\in\Nset, t\in\{t_p,\dots, t_{p+1}-1\}\}\ .
\end{equation}
Assume existence of a map $\alpha : \mathcal{T} \to [0,1]$ such that, for any $(p,t)\in\cT$, for any $k\in\cN$,
\begin{equation}
\label{exem}
k \in \cN_{p,t - t_p +1}
\Rightarrow
\sum_{l\in\cN_{p,t-t_p}} 
 \gamma_{k,l}(t) \geq\alpha(p,t)\ .
\end{equation}

Then the contraction rate of system \eqref{equb} as defined in Definition {\rm\ref{de15}} verifies:
\begin{equation}
\label{reca}
\rho\leq\limsup_{p\to +\infty}
\prod_{p'=1}^p
\left(
1-\prod_{t=t_{p'}}^{t_{p'+1}-1}\alpha(p',t)
\right)^{1/t_{p+1}}
\ .
\end{equation}
\end{theo}


Notice that, with the definition adopted in \eqref{cT}, there is indeed, for each $t\in\Nset$, a unique $p\in\Nset$ such that $(p,t)\in\cT$.

An important feature is that self-loops ($\gamma_{k,k}>0$) are not mandatory here, contrary to other previous contributions, see \cite{BER,MOR04a,BLO}.
This assumption is loosened up in \cite{TSI84,TSI86} for some, but not all, agents.
Example \ref{ex2} below presents an example where this is further weakened.
In particular, this feature permits to model leader/follower evolutions as well as leaderless networks within a unified framework.
 On this subject, see also Subsection \ref{se33} below.

Similarly, no positive uniform lower bound on the nonzero coefficients of $\Gamma(t)$ is required: requirement \eqref{exem} is sensibly weaker than the usual one in the literature, see \cite{CAO05,CAO06,OLS} and Example \ref{ex3} in the sequel.

\begin{proof}[Proof of Theorem \rm\ref{th1}]
One first states a monotonicity result for the diameter of the agent set along the solutions of \eqref{equb}.
\begin{lemm}
\label{le15}
For any trajectory of \eqref{equb}, one has, for any $t\in\Nset$,
\[
\Delta(t+1)\leq\Delta(t)\ .
\]
\end{lemm}
\begin{proof}
The proof of Lemma \ref{le15} comes from the fact that, the matrices $\Gamma(t)$ being stochastic, the map $t\mapsto\sup_{k\in\cN} x_k(t)$ (resp.\ $t\mapsto\inf_{k\in\cN} x_k(t)$) is non-increasing (resp.\ non-decreasing).
\end{proof}

One deduces directly from \eqref{recaa} that
\[
\Delta(t_{p+1})\leq \prod_{p'=1}^p \left(
1-\prod_{t_{p'}}^{t_{p'+1}-1} \alpha(p',t)
\right)\Delta(t_1)\ .
\]
Thus,
\[
\limsup_{p\to +\infty}
e^{\frac{1}{t_{p+1}}
\ln\left(
\frac{\Delta(t_{p+1})}{\Delta(0)}
\right)}
= \limsup_{p\to +\infty}
e^{\frac{1}{t_{p+1}}
\ln\left(
\frac{\Delta(t_1)}{\Delta(0)}
\right)}
e^{\frac{1}{t_{p+1}}
\ln\left(
\frac{\Delta(t_{p+1})}{\Delta(t_1)}
\right)}
\leq
\limsup_{p\to +\infty}
\prod_{p'=1}^p
\left (
1-\prod_{t=t_{  p'}}^{t_{  p'+1}-1}\alpha(p',t)
\right )^{1/t_{p+1}}\ .
\]
Clearly,
\[
\limsup_{p\to +\infty}
e^{\frac{1}{t_{p+1}}
\ln\left(
\frac{\Delta( t_{p+1}  )}{\Delta(0)}
\right)}
\leq
\limsup_{t\to +\infty} e^{\frac{1}{t} \ln \left ( \frac{\Delta(t)}{\Delta(0)} \right )}\ .
\]
Now, from the fact that $\Delta(t)$ is nonincreasing  and $t_{p+1}/t_{p+2} \rightarrow 1$ as $p \rightarrow + \infty$ one gets
\[
\limsup_{t\to +\infty} e^{\frac{1}{t} \ln \left ( \frac{\Delta(t)}{\Delta(0)} \right )}
\leq
\limsup_{p\to +\infty}
e^{\frac{1}{t_{p+2}}
\ln\left(
\frac{\Delta(t_{p+1})}{\Delta(0)}
\right)}
 = \limsup_{p \to + \infty} e^{\frac{1}{t_{p+1}}
\ln\left(
\frac{\Delta(t_{p+1})}{\Delta(0)}
\right)}
\]
(notice that the logarithmic expressions are not positive, due to the non-increasingness of $\Delta$ along time).
The conclusion is then immediate from the definition of $\rho$ given in Definition \ref{de15}.
\end{proof}

The next result is a specialisation of Theorem \ref{th1} for constant $\alpha$.

\begin{coro}
\label{co1}
Let the sequence of  communication graphs of system \eqref{equb} be $T$-sequentially connected.
Assume the existence of a constant map $\alpha$ in $[0,1]$ satisfying \eqref{exem}.
Then
\begin{equation}
\label{decr}
\rho\leq (1-\alpha^T)^{\frac{1}{T}}\ .
\end{equation}
\end{coro}
The previous results extends similar estimates found previously (see 
\cite{BER,CAO05,CAO06,NED}), as $\alpha$ does not have to bound from below the components of the matrices $\Gamma(t)$.

\begin{proof}[Proof of Corollary \rm\ref{co1}]

Assume without loss of generality $t_{p+1}-t_p \leq T$ for all $p \in \Nset$. Consequently $\limsup_{p \rightarrow + \infty} p T/ t_{p+1} \geq 1$.  Applying Theorem \ref{th1} with constant $\alpha$ yields for every $p\in\Nset$:
\begin{multline*}
\rho
\leq
\limsup_{p\to +\infty}
\prod_{p'=1}^p
\left(
1-\alpha^{t_{p'+1}-t_{p'}}
\right)^{1/t_{p+1}}
\leq
\limsup_{p\to +\infty}
\prod_{p'=1}^p
\left(\left(
1-\alpha^{T}
\right)^{1/T}
\right)^{T/t_{p+1}}\\
=
\limsup_{p\to +\infty}
\left(\left(
1-\alpha^T
\right)^{1/T}
\right)^{p T /t_{p+1}}
\leq (1-\alpha^T)^{1/T}\ .
\end{multline*}
where one has used the fact
that $T\mapsto (1 -\alpha^T)$ is increasing on $\Rset^+$, for any $\alpha\in [0,1]$.
Corollary \ref{co1} is thus proved.
\end{proof}


\begin{rema} A classical topic in linear algebra is the estimate of the second largest eigenvalue (in modulus) of a stochastic matrix for large dimensions.
In particular, as $n$ grows, Landau and Odlyzko \cite{landau} showed that the rate of convergence is of order $1 - 1/n^3$ (with $n$ being the number of agents)  for the equal neighbor time invariant model on undirected graph; see also results of the same nature in \cite{OLS}.
Our results can be applied to large systems as well. In particular, each given topology induces some kind of relation (typically an inequality) between tree-depth, weight of edges and number of agents.
This inequality can, in principle, be used to derive convergence rate estimates based on the number of agents. 
\end{rema}

We now provide several examples of systems with $n=3$ agents, in order to illustrate the two previous results.

\begin{exam}
\label{ex3}
As  a first example, consider the stationary system with $n=3$ agents given by:
\[
\Gamma= \Gamma(\varepsilon)\doteq\begin{pmatrix}
1/3 & 1/3 & 1/3\\
1/3 & 1/3 & 1/3\\
1/3 & 2/3-\varepsilon & \varepsilon
\end{pmatrix}\ ,
\]
for fixed $\varepsilon\in [0,1/3]$.
For $\varepsilon =1/3$, we obtain the equal neighbor averaging model corresponding to complete graph \cite{OLS}.
Spectral analysis argument shows that the actual value of the contraction rate $\rho$ is equal to $1/3-\varepsilon$.
Taking into account the fact that the coefficients are greater or equal to $\min\{1/3, 2/3-\varepsilon, \varepsilon\}=\varepsilon$, methods in \cite{CAO05,CAO06,OLS} yield an upper estimate of $\rho$ equal to $1-\varepsilon^2$, or even $1-\varepsilon\geq 2/3$ (taking into account the fact that the system under study is neighbor shared \cite{CAO05,CAO06} and adapting \cite[Lemma 2 and Theorem 1]{CAO05} to systems whose nonzero coefficients are at least $\varepsilon$).
 Now, Corollary \ref{co1} can be applied.
Indeed, the system appears as {\em 1-sequentially connected} --- as the first node participates with nonzero weight  to the evolution of all the agents ---, and one can take $\alpha=1/3$ as a lower-bound for these weights.
This gives an estimate of $\rho$ equal to $2/3$, which is better than the results obtained by the other methods.
\qed
\end{exam}

\begin{exam}
\label{ex1}
Consider a system with $n=3$ agents and dynamics defined by \eqref{equb} with
\begin{equation}
\label{Gamma}
\Gamma(t)
\doteq \left(
\begin{array}{c|c}
\alpha(t)M_1(t) & \begin{matrix}
\star & \star \\ \star & \star
\end{matrix}\\
\hline
\star &  \begin{matrix}
\star & \star
\end{matrix}
\end{array}
\right)
\text{ if $t\in 2\Nset$},\
\Gamma(t)
\doteq \left(
\begin{array}{c|c}
\alpha(t)M_2(t) & \begin{matrix}
\star \\ \star \\ \star
\end{matrix}
\end{array}
\right) \text{ if $t\in 2\Nset+1$}\ .
\end{equation}
Here, $\alpha\ :\ \Nset\to (0,1]$, $M_1\ :\ \Nset\to\Mset^{2, 1}$, $M_2\ :\ \Nset\to\Mset^{3, 2}$ (these sets have been defined in \eqref{Mpq}), and the stars stand for any nonnegative scalar numbers rendering the matrix $\Gamma(t)$ stochastic (for this to hold, the row-sums of $M_1(t)$, $M_2(t)$ have to be at most equal to $\alpha(t)^{-1}$).

 What is meant here is that for even $t$, $x_1(t)$ participates with a weight at least $\alpha(t)$ to the value of $x_1(t+1)$ and $x_2(t+1)$; and that for odd $t$, $x_1(t)$ and $x_2(t)$ participates with a {\em global} weight at least $\alpha(t)$ to the value of $x_1(t+1)$, $x_2(t+1)$ and $x_3(t+1)$.
This is precisely the assumptions needed to apply Theorem \ref{th1}, as detailed now.

Taking the first component $x_1$ as the root of the spanning-tree, one sees clearly that this system is 2-sequentially connected.
Application of Theorem \ref{th1} then yields the following estimate:
\begin{multline*}
\max_{t'=2t, 2t+1}\left\{
\max_{i=1,\dots, 3} x_i(t')-\min_{i=1,\dots, 3} x_i(t')
\right\}\\
\leq
\left(
1-\alpha(0)\alpha(1)
\right)\dots
\left(
1-\alpha(2t-2)\alpha(2t-1)
\right)
\left(
\max_{i=1,\dots, 3} x_i(0)-\min_{i=1,\dots, 3} x_i(0)
\right)\ ,
\end{multline*}
valid for any $t\in\Nset$.
When $\alpha$ is constant, Corollary \ref{co1} applies, and leads to:
\[
\max_{i=1,\dots, 3} x_i(t)-\min_{i=1,\dots, 3} x_i(t)
\leq (1-\alpha^2)^{\lfloor \frac{t}{2}\rfloor}
\left(
\max_{i=1,\dots, 3} x_i(0)-\min_{i=1,\dots, 3} x_i(0)
\right)\ .
\]

The following numerical experiment has been achieved.
A set of one thousand couples of stochastic matrices $\Gamma(1)$ and $\Gamma(2)$ are generated randomly (uniform law on $[0,1]$ is used for each coefficients, and the rows are afterward normalised), and the best estimates for $\alpha(1), \alpha(2)$ fulfilling the conditions above are then computed.
The actual contraction rate $\rho$ (which is the square-root of the maximal absolute value of the second largest eigenvalues $|\lambda_2(\Gamma(2)\Gamma(1))|$, see \cite[Proposition 1]{OLS}) is then compared to the upper bound $\trho$ deduced from Theorem \ref{th1} (that is $\sqrt{1-\alpha(1)\alpha(2)}$).
The corresponding histogram is represented in Figure \ref{fi1}.
\begin{figure}[h]
\begin{center}
\begin{tabular}{|c|c|c|c|c|c|c|}
\hline
0--5\% & 5--10\% & 10--15\% & 15--20\% & 20--25\% & 25--30\% & 30--35\% \\ 
1 & 21 & 86 & 161 & 212 & 197 & 132\\
\hline
35--40\% & 40--45\% & 45--50\% & 50--55\% & 55--60\% & 60--65\% & 65--70\% \\
88 & 52 & 30 & 13 & 4 & 2 & 1\\
\hline
\end{tabular}
\caption{Numerical test of Theorem \ref{th1}. Number of occurrences per value of the ratio $\rho/\trho$. See Example \ref{ex1} for details.}
\label{fi1}
\end{center}
\end{figure}
\qed
\end{exam}

Example \ref{ex1} shows that, although not tight, the bound may provide reasonable estimates.
Notice however that the previous comparison test is achieved only with 2-periodic systems (characterised by the second eigenvalue $\lambda_2(\Gamma(2)\Gamma(1))$), although Theorem \ref{th1} requires no specific assumption on the general time dependence.
An attempt to take into account the occurrence of several spanning-trees is proposed below (Section \ref{se55}).\\

\begin{exam}
\label{ex2}
We consider here a simple 2-periodic 3-agent system whose evaluation is not possible by the methods presented by previous works.
For $t\in\Nset$, we let
\[
\Gamma(2t)\doteq\begin{pmatrix}
1/2 & 1/2 & 0\\
1/2 & 1/2 & 0\\
0 & 0 & 1
\end{pmatrix},\quad
\Gamma(2t+1)\doteq\begin{pmatrix}
0 & 1/2 & 1/2\\
1/2 & 0 & 1/2\\
1 & 0 & 0
\end{pmatrix}\ .
\]
The matrices $\Gamma(2t+1)$ being deprived of any self-loop, the criteria from \cite{BER,MOR04a,BLO} cannot be applied.
Considering that the system is 2-sequentially connected (with $\cN^{1}_{p,0}=\{1\}$, $\cN^{1}_{p,1}=\{1,2\}$), use, as in Example \ref{ex1}, of Corollary \ref{co1} with $\alpha=1/2$ yields an estimate of the contraction rate as $(1-1/4)^{1/2}$, that is $\sqrt{3}/2\simeq 0.87$.
Indeed, the present example is an instance of Example \ref{ex1}, with $\alpha\equiv \frac{1}{2}$ in \eqref{Gamma}.
On the other hand, using as previously the second eigenvalue argument \cite{OLS}  of the product $\Gamma(2t+1)\Gamma(2t)$, the actual rate is found equal to $\sqrt{5/8}\simeq 0.79$.
 This value is smaller than $\sqrt{3}/2$, but it is computed under the restrictive hypothesis of periodicity.
\qed
\end{exam}

\begin{exam}
\label{ex35}
As a last illustration of Theorem \ref{th1}, an elementary time-varying 2-agent system is provided, for which no uniform-in-time lower bound on the nonzero coefficients of the state matrices exists.
This is a situation excluded from the previously published criteria.
Let
\[
\Gamma(2t)\doteq\begin{pmatrix}
1 & 0 & 0\\ 1 & 0 & 0\\ 0 & 1 & 0
\end{pmatrix},\
\Gamma(2t+1)\doteq\begin{pmatrix}
0 & 1 & 0\\ 1 & 0 & 0\\ \frac{1}{t} & 1-\frac{1}{t} & 0
\end{pmatrix}\ .
\]

This is clearly another special instance of Example \ref{ex1}.
Theorem \ref{th1} applies with $\alpha=1$, and yields a null contraction rate.
Indeed, finite-time convergence does occur, as
\[
\Gamma(2t)\Gamma(2t+1)=\begin{pmatrix}
1 & 0 & 0\\ 1 & 0 & 0\\ 1 & 0 & 0
\end{pmatrix}\ .
\]
\qed
\end{exam}


\subsection{Remarks on self-loops and delays}
\label{se33}

As noticed previously, self-loops are not assumed in the previous results.
However, it is a well-known fact that their absence may lead to non-convergence, as shown e.g.\ by the elementary example
\[
x_1(t+1)=x_2(t),\quad x_2(t+1)=x_1(t)
\]
whose solutions are either constant or oscillating, although the system is undoubtedly sequentially connected.

In fact, the assumptions of the main result of Section \ref{se4}, Theorem \ref{th1}, force the existence of a self-loop on the root of the spanning tree at each $t=t_p$, $p\in\Nset$.
This unique self-loop, together with the other hypotheses  (which impose rather information flux from upstream), turns out to be sufficient to enforce the convergence.
A similar remark will hold for the forthcoming results on systems with several spanning trees given below (see Section \ref{se55}, especially Theorem \ref{th3}).
This is a crucial point, as it drastically conditions the search for spanning-trees.\\

On the other hand, it is quite evident that the information transfer between the agents may be subject to delays.
This feature does not present specific difficulties a priori, because it can be treated the same way, via an augmentation of the state vector.
Some past state values are then included in the definition of the diameter and of the contraction rate which are considered (see Definitions \ref{de1} and \ref{de15}), but this has essentially no consequence on the meaning of this  latter quantity.

Generally speaking, delays cannot suppress sequential connectivity, except if they concern the unique mandatory self-loop, located at the root at the initial time instant (see above).
On the other hand, they may change the values of the weights $\alpha(t)$ and thus modify the decay estimates.
Also, it is rather likely that the decay rate estimate are nondecreasing with respect to any delay.

A more precise study of the quantitative influence of the delays on the convergence speed could be tackled by similar tools, 
this feature is out of the scope of the present paper.
For simplicity, we limit ourselves to a simple example, for which analytical results are easily computed.

\begin{exam}
Consider the four following systems:
\begin{subequations}
\begin{gather}
\label{sima}
x_1(t+1) = \frac{1}{2} x_1(t)+\frac{1}{2} x_2(t),\quad
x_2(t+1) = \frac{1}{2} x_1(t)+\frac{1}{2} x_2(t)\ ;\\
\label{simb}
x_1(t+1) = \frac{1}{2} x_1(t)+\frac{1}{2} x_2(t),\quad
x_2(t+1) = \frac{1}{2} x_1(t-1)+\frac{1}{2} x_2(t-1)\ ;\\
\label{simc}
x_1(t+1) = \frac{1}{4} x_1(t)+\frac{1}{4} x_1(t-1)+\frac{1}{2} x_2(t), \quad
x_2(t+1) = \frac{1}{2} x_1(t-1)+\frac{1}{2} x_2(t)\ ;\\
\label{simd}
x_1(t+1) = \frac{1}{2} x_1(t-1)+\frac{1}{2} x_2(t),\quad
x_2(t+1) = \frac{1}{2} x_1(t)+\frac{1}{2} x_2(t-1)\ .
\end{gather}
\end{subequations}

The delay-free system \eqref{sima} is 1-sequentially connected, and the analysis conducted above yields the estimate $\trho=\frac{1}{2}\geq \rho=0$.

The three remaining systems possess delayed terms $x_1(t-1)$ and $x_2(t-1)$.
Introducing $x_3(t)\doteq x_1(t-1)$, $x_4(t)\doteq x_2(t-1)$, they can be written $x(t+1)= \Gamma x(t)$ where $x\doteq (x_1,x_2,x_3,x_4)\t$ and
\[
\Gamma
\doteq
\begin{pmatrix}
\frac{1}{2} & \frac{1}{2} & 0 & 0\\
0 & 0 & \frac{1}{2} & \frac{1}{2}\\
1 & 0 & 0 & 0\\
0 & 1 & 0 & 0
\end{pmatrix},\quad
\Gamma
\doteq
\begin{pmatrix}
\frac{1}{4} & 0 & \frac{1}{4} & \frac{1}{2}\\
\frac{1}{2} & 0 & 0 & \frac{1}{2}\\
1 & 0 & 0 & 0\\
0 & 1 & 0 & 0
\end{pmatrix},\quad
\Gamma
\doteq
\begin{pmatrix}
0 & \frac{1}{2} & \frac{1}{2} & 0 \\
\frac{1}{2} & 0 & 0 & \frac{1}{2}\\
1 & 0 & 0 & 0\\
0 & 1 & 0 & 0
\end{pmatrix}
\]
respectively for \eqref{simb}, \eqref{simc} and \eqref{simd}.
System \eqref{simb} is 3-sequentially connected with $\alpha(0)=\alpha(1)=\alpha(2)=\frac{1}{2}$, so that $\trho=\frac{\sqrt[3]{7}}{2}\simeq 0.96 \geq \rho=\frac{1}{2}$.
System \eqref{simc} is 2-sequentially connected with $\alpha(0)=\frac{1}{4}$ and $\alpha(1)=\frac{1}{2}$, and $\trho=\frac{\sqrt{7}}{2\sqrt{2}}\simeq 0.94 \geq \rho=\frac{1}{2}$.
Last, system \eqref{simd} is not sequentially connected, due to the absence of self-loop (zero diagonal).
This is corroborated by the fact that $\rho=1$, so convergence does not occur.

As can be seen, the index of sequential connectivity is not systematically the sum of this index for the delay-free case (1 here) and the sum of the values of the delays.
\end{exam}

\section{Communication graphs spanned by several spanning-trees} 
\label{se55}

\subsection{Sequential connectivity with several spanning-trees}
\label{se41}

When several spanning-trees emerge in the communication graph (either simultaneously, or successively), the previous analysis may happen to be conservative.
We now face the issue of how to tackle this feature.

An extension of the notion of sequential connectivity introduced in Section \ref{se3} is first constructed, 
analogously to Definitions \ref{de4bis} and \ref{de5}.

\begin{defi}[Sequential connectivity of finite graph sequences by multiple spanning-trees]
\label{de4bismul}
A finite sequence of $T$ graphs with common nodes $\cA_1, \cA_2, \ldots, \cA_T$ is said to be {\em sequentially connected 
by $m$ spanning-trees} ($m\in\Nset$) if
there exist  nodes $k^1, k^2, \ldots, k^m \in \cN$ and iterations given by:
\begin{eqnarray} 
\label{spantreemultiple}
\cN^j_0 &=& \{ k^j \} \nonumber \\
\cN^j_{t} & \subseteq & \cN^j_{t-1}  \cup \neig ( \cN^j_{t-1}, \cA_{t} ) \quad t = 1, \ldots, T \nonumber
\end{eqnarray}
satisfy $\cN^j_{T}= \cN$ for all $j \in \cI\doteq\{ 1, 2, \ldots , m \}$. \\
\end{defi}

Similarly we define the property for infinite graph sequences:
\begin{defi}[$T$-sequential connectivity by multiple spanning-trees]
\label{de10}
An infinite sequence of graphs $\cA (t)$, $t \in \Nset$, is said {\em $T$-sequentially connected by $m$ spanning-trees} if there 
exists a strictly increasing integer sequence $t_p$, $p\in\Nset$, fulfilling \eqref{defisup} and such that for all
$p \in \Nset$, each graph sub-sequence
\[
\cA ( t_p ),  \ldots, \cA (t_{p+1}-1) \]
is sequentially connected by $m$ spanning-trees.
\end{defi}

The following result extends Lemma \ref{le0}.
As the latter, it is directly deduced from Lemma \ref{le1}, so detailed proof is omitted.
\begin{lemm}
\label{le3}
Let the finite sequence of communication graphs $\cA(0),\dots, \cA(T-1)$ of system \eqref{equb} be sequentially connected by $m$ spanning-trees, and let $\cN^j_0,\dots,\cN^j_T$ be the sets corresponding to the $j$-th spanning-tree (see Definition {\rm\ref{de4bismul}}).
For each $k\in\cN$, consider sets $\cN^{k,j}_t\subseteq\cN^j_t$ such that $j\neq j'\Rightarrow\cN^{k,j}_t\cap\cN^{k,j'}_t=\emptyset$.
Assume the existence of maps $\alpha_{i,j}: [0,T-1]\to [0,1]$ such that, for any $t=0,\dots, T-1$, any  $i,j\in\cI$ and any $k\in\cN$,
\begin{equation}
\label{exemc}
k \in \cN^i_{t+1} 
\Rightarrow
\sum_{l\in\cN^{k,j}_{t} } 
 \gamma_{k,l}(t) \geq\alpha_{i,j}(t)\ ,
\end{equation}
for given maps $\alpha_{i,j}: [0,T-1]\to [0,1]$.
Then $A(t)\doteq (\alpha_{i,j}(t))_{(i,j)\in\cI\times\cI}$ is a sub-stochastic matrix in $\Rset^{m\times m}$ and
\begin{equation}
\label{recc}
\Delta(T) \leq 
\left( 1-\max_{j\in\cI}\bfo_j\t A(T-1)\dots A(0)\bfo
\right) \, \Delta(0)\ ,
\end{equation}
where $\bfo$ and $\bfo_j$ have been defined in  the Notation section.
\end{lemm}

The sets $\cN^j_t$ generalise the notion introduced in Section \ref{se4}: here, each set $\cN^j_t$ is constituted by agents attained at most at time $t$ by the $j$-th spanning-tree.
Assumption \eqref{exemc} fixes a lower bound $\alpha_{i,j}$ to the total weight applied by each agent in $\cN^i_{t+1}$ to agents in $\cN^j_t$.

When an agent is member, for given $t$, of sets $\cN^j_t$ for more than one value of $j$, it is necessary to decide, in the construction of contraction estimates relative to the $x_k$ update equation, to which of them its influence is attributed; this choice could actually vary according to the considered update equation. This is the reason why subsets $\cN^{k,j}_t$ disjoint for different $j$s are introduced.

\begin{rema}
When in the statement of Lemma \ref{le3} the sequential spanning-trees corresponding to two distinct values of $i, i'$ are identical, then the scalar quantities $\bfo_i\t A(T-1)\dots A(0)\bfo$ and $\bfo_{i'}\t A(T-1)\dots A(0)\bfo$ are equal --- at least if the $\alpha_{i,j}(t)$ are chosen identical for all $j\in\cI$.

On the other hand, different choices in the attribution of arcs to one or another of the $m$ developing spanning-trees (that is on the definition of the sets $\cN^{k,j}_t$) may lead to different choices for these coefficients, and consequently to different estimates.
In this respect, adding virtual sequential spanning-trees may allow to improve the convergence speed estimate, see Example \ref{ex5} below.
\qed
\end{rema}

\begin{rema}
As in Remark \ref{re2} for Lemma \ref{le0}, one shows easily that, under the hypotheses of Lemma \ref{le3}, one has similarly, for any $j\in\cI$,
\[
\forall l\in\cN,\ \frac{\partial x_l(T)}{\partial x_k(0)} \geq \bfo_j\t A(T-1)\dots A(0)\bfo_k\ ,
\]
where the index $k$ represents any index of the roots of the $m$ spanning-trees.
Indeed, for any $t=0,\dots, T-1$ and any $j\in\cI$,
\[
\forall l\in\cN^j_{t+1},\ \frac{\partial x_l(t+1)}{\partial x_k(0)} \geq \bfo_j\t A(T-1)\dots A(0)\bfo_k\ ,
\]
which gives the previous estimate when $t=T-1$, as $\cN^j_T=\cN$ for all $j\in\cI$.
\qed
\end{rema}

\subsection{Results on contraction rate estimate}
\label{se42}

We now come to the key result of Section \ref{se55}, which is also the most powerful of the paper.
\begin{theo}
\label{th3}
Let the sequence of  communication graphs of system \eqref{equb} be $T$-sequentially connected by $m$ spanning-trees.
Accordingly, denote $t_p$ the corresponding increasing  integer sequence of spanning-tree completion (see Definition {\rm\ref{de10}}); $\cN_{p,t-t_p}^j$, $t=t_p,\dots, t_{p+1}$, the sets corresponding to the $m$ spanning-trees connecting sequentially the graph sub-sequences $\cA(t_p),\dots, \cA(t_{p+1}-1)$, $p\in\Nset$ (see Definition {\rm\ref{de4bismul}}); and $\cT$ the set defined in \eqref{cT}.
Furthermore, for each $k\in\cN$, consider sets $\cN^{k,j}_{p,t-t_p}\subseteq\cN^j_{p,t-t_p}$ such that $j\neq j'\Rightarrow \cN^{k,j}_{p,t-t_p}\cap\cN^{k,j'}_{p,t-t_p} =\emptyset$, and assume the existence of maps $\alpha_{i,j}: \cT\to [0,1]$, $(i,j)\in\cI\times\cI$, such that, for any $(p,t)\in\cT$, any  $i,j\in\cI$ and any $k\in\cN$,
\begin{equation}
\label{exemcc}
k \in \cN^i_{p,t-t_p+1} 
\Rightarrow
\sum_{l\in\cN^{k,j}_{p,t-t_p}} 
 \gamma_{k,l}(t) \geq\alpha_{i,j}(p,t)\ .
\end{equation}

Then $A(p,t)\doteq (\alpha_{i,j}(p,t))_{(i,j)\in\cI\times\cI}$ is a sub-stochastic matrix in $\Rset^{m\times m}$ and the contraction rate of system \eqref{equb} as defined in Definition {\rm\ref{de15}} verifies:
\begin{equation}
\label{reccc}
\rho\leq\limsup_{p \to +\infty}
\prod_{p'=1}^p
\left(
1-\max_{j\in\cI} \bfo_j\t A(p',t_{p'+1}-1)\dots A(p',t_{p'})\bfo
\right)^{1/t_{p+1}}\ .
\end{equation}
\end{theo}


\begin{proof}[Proof of Theorem {\rm\ref{th3}}]
Due to the fact that the sets $\cN_{k,j}(t)$ are pairwise disjoint for different values of $j$, one has, for any $i\in\cI$, any $t\in\Nset$,
\[
\sum_{j\in\cI}\alpha_{i,j}(p,t-t_p)
\leq \sum_{j\in\cI}\sum_{l\in\cN_{k,j}(t)}\gamma_{k,l}(t)
\leq \sum_{l\in\cN}\gamma_{k,l}(t)=1\ ,
\]
for any $k\in\cN^i_{p,t-t_p+1}$.
This proves the first part of the statement.

As in the proofs of Theorem \ref{th1} above, it suffices essentially to establish \eqref{recc} when $t$ is a multiple of $T$.
Applying Lemma \ref{le1} on the integer interval $[t_p, t_{p+1}]$ with $\cM_i(t)=\cN^i_{p,t-t_p}$, $\cM_{k,j}(t)=\cN^{k;j}_{p,t-t_p}$, $c_{i,j}(t)=\alpha_{i,j}(p,t)$ yields:
\[
\Delta_i( t_{p+1})
\leq
\bfo_i\t \left(
 A(p,t_{p+1}-1)\dots A(p,t_p) \Delta_\cN(t_p)
+\left(
\vphantom{\int}
\bfo
-  A(p,t_{p+1}-1)\dots A(p,t_p) 
\bfo
\right)
\Delta(t_p)\right)\ .
\]
Here,  
the definition of $\Delta_\cN$ depends upon $p$ and is as follows:
\[
\Delta_\cN(t)\doteq\begin{pmatrix}
\Delta_1(t)\\ \vdots \\ \Delta_m(t)
\end{pmatrix}, \quad 
\Delta_i(t)\doteq\max_{k\in  \cN^i_{p,t-t_p}} x_k(t) -\min_{k\in \cN^i_{p,t-t_p}} x_k(t),\quad
 t=t_p,\dots,t_{p+1}\ .
\]
By assumption, the existence of the $m$ spanning-trees means that:
\[
\Delta_i( t_p)=0 \text{ and } \Delta_i( t_{p+1})=\Delta( t_{p+1}),\ i\in\cI\ .
\]
One thus deduces that, for all $i\in\cI$,
\[
\Delta( t_{p+1})
=\Delta_i( t_{p+1})
\leq
\bfo_i\t \left(
\vphantom{\int}
\bfo
-  A(p,t_{p+1}-1)\dots A(p,t_p)\bfo
\right)
\Delta( t_p)\ .
\]
Thus,
\[
\Delta( t_{p+1}) \leq \left(
1- \bfo_i\t  A(p,t_{p+1}-1)\dots A(p,t_p)\bfo
\right)
\Delta( t_p)\ .
\]
 The proof is then achieved as for Theorem \ref{th1}.
\end{proof}

\begin{exam}
\label{ex5}
We come back to the analysis of Example \ref{ex3}, now with the help of Theorem \ref{th3}.
One may distinguish three spanning-trees occurring on each time interval of unit length (in other words, the system is 1-sequentially connected by 3 spanning-trees), with root at each of the agents.
With this point of view, $\cI=\cN= \{ 1, 2, 3 \}$  and  $t_p=p$.
With the notation of Theorem \ref{th3}, one may put:
\[
\cN^j_{p,0}=\{j\},\ \cN^j_{p,1} =\cN=\{1,2,3\}, \text{ for } j=1,2, 3\ . 
\]
This is the simplest case, where the sets $\cN^j_{p,0}$ are pairwise disjoint
, so one takes
\[
\cN^{k,j}_{p,0}\doteq\cN^j_{p,0},\ j=1,2,3\ .
\]

We now form the functions $\alpha_{i,j}$ as defined in the statement of Theorem \ref{th3}, and the corresponding matrix $A$.
By definition, one should have for any $i,j\in\{1,2,3\}$ (see \eqref{exemcc}):
\[
\forall k\in\{1,2,3\},\forall p\in\Nset,\ \gamma_{k,j}\geq\alpha_{i,j}(p,p)\ ,
\]
where $\Gamma= (\gamma_{i,j})_{(i,j)\in\cI\times\cI}$ is given in Example \ref{ex3} above.
One thus takes
\[
\alpha_{i,j}(p,p)\doteq\min_{k=1,2,3}\gamma_{k,j},\ i=1,2,3\ ,
\]
that is: $\alpha_{i,1}=1/3$, $\alpha_{i,2}=\min\{1/3, 2/3-\varepsilon\}$, $\alpha_{i,3}=\min\{1/3,\varepsilon\}$, or again
\[
A\doteq\begin{pmatrix}
1/3 & \min\{1/3, 2/3-\varepsilon\} & \min\{1/3,\varepsilon\}\\
1/3 & \min\{1/3, 2/3-\varepsilon\} & \min\{1/3,\varepsilon\}\\
1/3 & \min\{1/3, 2/3-\varepsilon\} & \min\{1/3,\varepsilon\}
\end{pmatrix}\ .
\]

Applying then formula \eqref{recc} leads to an estimate of the actual contraction rate equal to
\[
1-\left(
\vphantom{\sum}
1/3+\min\{1/3, 2/3-\varepsilon\}+\min\{1/3,\varepsilon\}
\right)
= 1/3-\varepsilon\ .
\]
In this example, the method ensuing from Theorem \ref{th3} thus generates the exact value of the contraction rate $\rho$.

Considering now only the two first spanning-trees (with $\cN^j_{p,0}=\{j \}$, $\cN^j_{p,1}=\cN=\{1,2,3\}$ for
all $p \in \Nset$ and all $j \in \{1,2\}$; then $\alpha_{i,1}=1/3$, $\alpha_{i,2}=\min\{1/3, 2/3-\varepsilon\}=1/3$, $i=1,2$) gives a worse estimate, namely $1/3$. 
Similarly, considering the first and third, or the second and third, spanning-trees yields $2/3-\varepsilon$.
These estimates are different, tighter than 2/3, the value obtained in Example \ref{ex3} when considering a unique spanning-tree, but not optimal.
\qed
\end{exam}

We refine further in the sequel the analysis of systems spanned by several spanning-trees, and examine respectively in Sections \ref{se6} and \ref{se7} the cases of spanning-trees propagating consecutively and simultaneously.

\subsection{Application to systems with successive spanning-trees}
\label{se6}

We now consider the case where the several emerging spanning-trees have common root and possess certain order property.
We mean by this that the dates at which each spanning-tree reach an agent are  interlaced {\em independently\/} from the agent.
Otherwise said, the ``wavefronts" corresponding to each spanning-tree spread in a concentrical manner.
Up to renaming, one may label 1 the first spanning-tree, 2 the next one and so on\ \dots, and the order property simply reads (reasoning on each interval $[t_p, t_{p+1}]$, we omit the index $p$):
\[
\forall j,j'\in\cI, \forall t=0,\dots, T-1,\quad j\leq j'
\Rightarrow
\cN^{j'}_t\subseteq\cN^j_t\ ,
\]
and thus, by construction, the following inequalities hold, for any $t\in[t_p,t_{p+1}]$:
\[
\Delta_m(t)\leq\dots\leq\Delta_j(t)\leq\dots\leq\Delta_1(t)\ .
\]


It is thus systematically more fruitful to attribute any contribution in the right-hand side of \eqref{equb} to the set $\cN_i$ with largest index $i$ to which it belongs  --- because the corresponding estimate is tighter.
In particular, it is beneficial to choose $\alpha_{i,j}\equiv 0$ for $i<j$, thus leading to lower-triangular matrices $A$ in Theorem \ref{th3}.

We provide now an illustration of this configuration.

\begin{exam}
\label{ex6}
For a fixed scalar $\gamma\in [0,1]$, consider the time-invariant system of $n$ agents described by
\[
x_1(t+1)=x_1(t),\quad x_i(t+1)=\gamma x_{i-1}(t)+ (1-\gamma)x_i(t),\ i=2,\dots, n\ .
\]
The corresponding matrix $\Gamma$ is lower-triangular and admits, apart from 1, a unique eigenvalue, namely $1-\gamma$, with degree $n-1$.
The actual value of the contraction rate is thus $\rho=1-\gamma$.

For any positive integer $q$, one may consider that the communication graph is spanned by $q$ distinct spanning-trees, departing from agent 1 at time 0, then 1, 2 and so on, up to $q-1$, and attaining agent $n$ at time $n-1$, $n$, up to $n+q-2$.
The duration of this process is thus $T\doteq n+q-2$, and the system may be seen as ``$T$-sequentially connected by $q$ (distinct) spanning-trees".
Coherently with the previous notations, we let $t_p=pT$ and  consider the sets $\cN^i_{p,t}$, $i\in\cI\doteq\{1,\dots, q\}$, defined by:
\[
\cN^1_{p,t}=
\begin{cases}
\{1\} & \text{ for } t\leq 0,\\
\{1,\dots,t+1\} & \text{ for } t=0,\dots, n-1,\\
\{1,\dots,n\}=\cN & \text{ for } t=n-1,\dots, n+q-2,
\end{cases}
\text{ and }\quad
\cN^{i+1}_{p,t+1}=\cN^i_{p,t} \text{ for } i=2,\dots, q-1\ .
\]
Following the progression of each spanning-tree, one shows that one may take for $A$ (in $\Rset^{q\times q}$) the formulas depicted in Figure \ref{fi35} (see also in Appendix the details of the proof of Lemma \ref{le4} below).

\begin{figure}
\small
\begin{gather*}
A(0)=\begin{pmatrix}
\gamma & 0 & \dots &0\\
0 & 1 & \dots & 0\\
\vdots & \vdots& & \vdots \\
0 & 0 & \dots & 1
\end{pmatrix},\
A(1)=\begin{pmatrix}
\gamma & 0 & \dots & 0\\
1-\gamma & \gamma & \dots & 0\\
\vdots & \vdots & & \vdots \\
0 & 0 & \dots & 1
\end{pmatrix},\dots,\
A(n-2)=\begin{pmatrix}
\gamma & 0 & \dots & 0\\
1-\gamma & \gamma & \dots & 0\\
0 & 1-\gamma & \dots & 0 \\
\vdots & & & \vdots\\
0 & 0 & \dots & 1
\end{pmatrix}\ ,\\
\hspace{-.3cm}
A(n-1)=\begin{pmatrix}
1 & 0 & \dots & 0\\
1-\gamma & \gamma & \dots & 0\\
0 & 1-\gamma & \dots & 0 \\
\vdots & & & \vdots\\
0 & 0 & \dots & 1
\end{pmatrix},\
A(n)=\begin{pmatrix}
1 & 0 & \dots & 0\\
0 & 1 & \dots & 0\\
0 & 1-\gamma & \dots & 0 \\
\vdots & & & \vdots\\
0 & 0 & \dots & 1
\end{pmatrix},\dots,\
A(n+q-3)=\begin{pmatrix}
1 & 0 &  \dots\\
0 & 1 & \dots\\
\vdots & & &\\
0 &  \dots & 1 & 0 \\
0 &  \dots & 1-\gamma & \gamma
\end{pmatrix}\ .
\end{gather*}
\caption{Matrices $A$ obtained in Example {\rm\ref{ex6}} (case $n<q$).}
\label{fi35}
\end{figure}

Let us explain these formulas.
From $t=0$ to $t=n-1$, the first spanning-tree spreads from agent 1 to agent $n$;
for the elements of the latter, the right-hand side of the state equation is composed by element already touched by the information flow (with coefficient $\gamma$) and some newly touched element, which consequently does not contribute to the right-hand side. We may therefore choose, $\cN^{1,1}_{p,t} = \cN^1_{p,t}$ and
$\cN^{1,j}_{p,t}= \emptyset$ for $j=2 \ldots q$. 
This gives rise to the identities $\alpha_{1,1}(t)=\gamma$ and $\alpha_{1,j}(t)=0$ for $j\in\cI\setminus\{1\}$ for $t=0,\dots, n-2$.
At time $t=n-1$, one has $\cN^1_{p,t}=\cN$ and the expansion of this set is completed, so {\em all\/} the terms in the right-hand side come from inside $\cN^1_{p,t}$.
Thus, again by letting $\cN^{1,1}_{p,t}=\cN$ and $ \cN^{1,j}_{p,t}= \emptyset$  for $t=n-1,\dots, n+q-3$ and $j= 2 \ldots q$, one has $\alpha_{1,1}(t)=1$ and $\alpha_{1,j}(t)=0$.

The second spanning-tree departs from the root at $t=1$, therefore letting $\cN^{2,2}_{p,0}=\cN^{2}_{p,0}$
and $\cN^{2,j}_{p,0} = \emptyset$ yields $\alpha_{2,2}(t)=1$ and $\alpha_{2,j}(t)=0$ for $j\in\cI\setminus\{2\}$ for $t=0$.
Then at $t=1$, $\cN^2_{p,t}=\{1\}$ and $\cN^2_{p,t+1}=\{1,2\} = \cN^2_{p,t}\cup\left(
\cN^1_{p,t}\setminus\cN^2_{p,t}
\right)$.
More precisely, the corresponding right-hand side comprises two terms as before: a contribution, with coefficient $\gamma$, due to agents already attained by the second spanning-tree, plus a term, with coefficient $1-\gamma$, due to a term coming from an agent not yet touched by the second tree, {\em but already by the first one}.
We let $\cN^{2,1}_{p,1}=\{ 2 \}$, $\cN^{2,2}_{p,1}= \{1\}$, $\cN^{2,j}_{p,1}= \emptyset$ for $j=3 \ldots q$;
this explains that for $t=1$ one has: $\alpha_{2,1}(t)=1-\gamma$, $\alpha_{2,2}(t)=\gamma$ and $\alpha_{2,j}(t)=0$ for $j\in\cI\setminus\{1,2\}$.
Similarly, we define $\cN^{2,1}_{p,t}=\{t+1\}$, $\cN^{2,2}_{p,t} = \{ 1 \ldots t \}$ and 
$\cN^{2,j}_{p,t}= \emptyset$ for $t=2 \ldots n$ where the second spanning-tree in turn is completed.
Again one obtains $\alpha_{2,1}(t)=1-\gamma$, $\alpha_{2,2}(t)=\gamma$ and $\alpha_{2,j}(t)=0$ for $j\in\cI\setminus\{1,2\}$. Then for subsequent $t$'s, we let $\cN^{2,2}_{p,t}=\cN$ and $\cN^{2,j}=
\emptyset$ for $j \neq 2$. So that indeed $\alpha_{2,2}(t) = 1$ and $\alpha_{2,j}(t)=0$ for $j\neq 2$.

Last, the other spanning-trees appear one by one, and share with their predecessor the same relation than the second one with the first one.
This explains the formulas given, until completion of the $q$-th one, at time $t=T$.
The analysis conducted above leads overall to the matrices shown in Figure \ref{fi35}, which corresponds to the case $n<q$ (the first spanning-tree is completed at $t=n$, before the departure of the $q$-th spanning-tree, at $t=q$).
The case $n\geq q$ is similar.

For the case of $n=3$ agents, formula \eqref{recc} in Theorem \ref{th3} then yields the following estimates, denoted $\trho_q$:

\begin{itemize}

\item[$\bullet$]
for $q=1$ (corresponding to the method of Theorem \ref{th1}):
\[
A(0)=A(1)=\gamma\ ,
\]
so $\trho_1=\sqrt{1-A(1)A(0)}=\sqrt{1-\gamma^2}$.

\item[$\bullet$]
for $q=2$:
\[
A(0)=\begin{pmatrix}
\gamma & 0\\
0 & 1
\end{pmatrix},\
A(1)=\begin{pmatrix}
\gamma & 0\\
1-\gamma & \gamma
\end{pmatrix},\
A(2)=\begin{pmatrix}
1 & 0\\
1-\gamma & \gamma
\end{pmatrix}\ ,
\]
and
\[
\trho_2=\left(1-\max_{i=1,2} \bfo_i\t A(2)A(1)A(0)\bfo \right)^{1/3}
= \left(
1-\max\{\gamma^2;3\gamma^2-2\gamma^3\}
\right)^{1/3}
= \left(
1-3\gamma^2+2\gamma^3
\right)^{1/3} \ .
\]

\item[$\bullet$]
for $q=3$:
\[
A(0)=\begin{pmatrix}
\gamma & 0 & 0\\
0 & 1 & 0\\
0 & 0 & 1
\end{pmatrix},\
A(1)=\begin{pmatrix}
\gamma & 0 & 0\\
1-\gamma & \gamma & 0\\
0 & 0 & 1
\end{pmatrix},\
A(2)=\begin{pmatrix}
1 & 0 & 0\\
1-\gamma & \gamma & 0\\
0 & 1-\gamma & \gamma
\end{pmatrix},\
A(3)=\begin{pmatrix}
1 & 0 & 0\\
0 & 1 & 0\\
0 & 1-\gamma & \gamma
\end{pmatrix}\ ,
\]
whence:
\begin{eqnarray*}
\trho_3
& = &
\left(1-\max_{i=1,2,3} \bfo_i\t A(3)A(2)A(1)A(0)\bfo \right)^{1/4}\\
& = &
\left(
1-\max\{\gamma^2;3\gamma^2-2\gamma^3;3\gamma^4-8\gamma^3+6\gamma^2\}
\right)^{1/4}
= \left(
1-3\gamma^4+8\gamma^3-6\gamma^2
\right)^{1/4}\ .
\end{eqnarray*}
\end{itemize}
 
The values obtained approximates the exact value $1-\gamma$ with increasing precision, as seen in Figure \ref{fi3}.
\begin{figure}
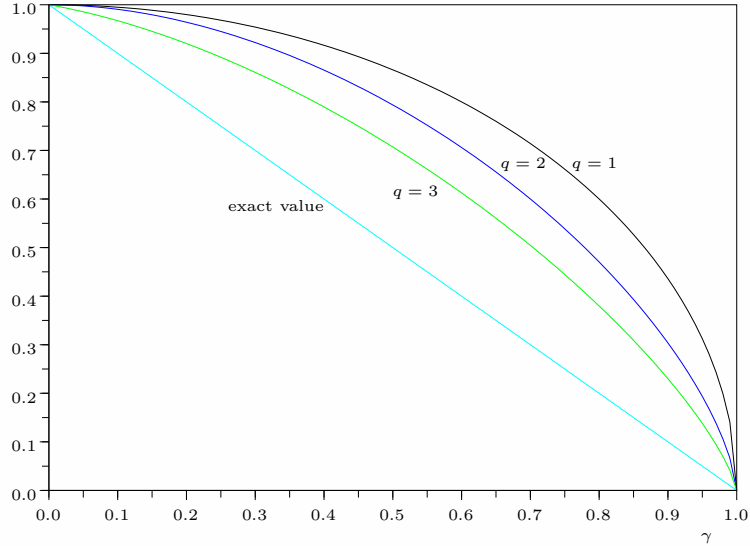

\begin{center}
\input curves1.pstex_t
\end{center} 
\caption{Approximations of the contraction rate as functions of $\gamma$, for different uses of Theorem \ref{th3}.
See text of Example \ref{ex6} for details.}
\label{fi3}
\end{figure}
These successive improvements are of course consequence of a richer and richer analysis, including more and more settling spanning-trees.

The question of the limiting behaviour when $q$ goes to infinity is of course intriguing: is the exact value found asymptotically?
It turns out that the answer is positive, as stated now in the general case of a system with $n$ agents.

\begin{lemm}
\label{le4}
The value of $\trho_q$ is given by the following formula:
\[
\trho_q^{n+q-2}
=
1-\gamma^{n-1}
\sum_{i=0}^{q-1}
(1-\gamma)^i\begin{pmatrix}
n+i-2 \\ n-2
\end{pmatrix}
= \frac{1}{(n-2)!} \gamma^{n-1} \frac{d^{n-2}}{d\delta^{n-2}}
\left. \left[
\frac{\delta^{n+q-2}}{1-\delta}
\right] \right|_{\delta=1-\gamma}\ .
\]
Consequently, $\trho_q$ tends towards $\rho= 1-\gamma$ when $q\to +\infty$, and more precisely
\[
\trho_q = \rho+(n-2)(1-\gamma) \frac{\ln q}{q} + o\left(
\frac{\ln q}{q}
\right)\ .
\]
\end{lemm}

A proof of Lemma \ref{le4} is presented in Appendix.
The calculations have been checked independently by the authors, using symbolic computation tool.

Although presently limited to special class of examples, Lemma \ref{le4} is rather promising: it establishes that tight estimates may be accessed to, when employing large number of settling spanning-trees in the analysis.
Extensions are in progress to cover more general cases.
\qed
\end{exam}

\subsection{Application to systems with concomitant spanning-trees}
\label{se7}

The example previously shown exploit drastically the fact that the different spanning-trees occur one after another.
We show here that, otherwise, the techniques of Theorem \ref{th3} may provide deceivingly weak results.

\begin{exam}
\label{ex55}
To illustrate this, we consider a system with $n=6$ agents, $T$-sequentially connected for $T=5$.
For fixed $\gamma\in (0,1/2)$, the latter is defined by taking stochastic matrices such that:
\begin{subequations}
\label{G55}
\begin{gather}
\Gamma(pT)
\geq\begin{pmatrix}
\gamma & 0 & 0 & 0 & 0 & 0\\
\gamma & 0 & 0 & 0 & 0 & 0\\
0 & 0 & 0 & 0 & 0 & 0\\
0 & 0 & 0 & 0 & 0 & 0\\
0 & 0 & 0 & 0 & 0 & 0\\
0 & 0 & 0 & 0 & 0 & 0
\end{pmatrix},\
\Gamma(pT+1)
\geq\begin{pmatrix}
\gamma & 0 & 0 & 0 & 0 & 0\\
0 & \gamma & 0 & 0 & 0 & 0\\
0 & \gamma & 0 & 0 & 0 & 0\\
0 & \gamma & 0 & 0 & 0 & 0\\
0 & 0 & 0 & 0 & 0 & 0\\
0 & 0 & 0 & 0 & 0 & 0
\end{pmatrix},\
\Gamma(pT+2)
\geq\begin{pmatrix}
\gamma & 0 & 0 & 0 & 0 & 0\\
0 & \gamma & 0 & 0 & 0 & 0\\
0 & 0 & \gamma & 0 & 0 & 0\\
0 & 0 & 0 & \gamma & 0 & 0\\
0 & 0 & \gamma & 0 & 0 & 0\\
0 & 0 & 0 & \gamma & 0 & 0
\end{pmatrix},\\
\Gamma(pT+3)
\geq\begin{pmatrix}
\gamma & 0 & 0 & 0 & 0 & 0\\
0 & \gamma & 0 & 0 & 0 & 0\\
0 & 0 & \gamma & 0 & 0 & 0\\
0 & 0 & 0 & \gamma & 0 & 0\\
0 & 0 & 0 & 0 & \gamma & \gamma\\
0 & 0 & 0 & 0 & \gamma & \gamma
\end{pmatrix},\
\Gamma(pT+4)
\geq\begin{pmatrix}
\gamma & 0 & 0 & 0 & 0 & 0\\
0 & \gamma & 0 & 0 & 0 & 0\\
0 & 0 & 0 & 0 & \gamma & 0\\
0 & 0 & 0 & 0 & 0 & \gamma\\
0 & 0 & 0 & 0 & \gamma & 0\\
0 & 0 & 0 & 0 & 0 & \gamma
\end{pmatrix}
\end{gather}
\end{subequations}
for all $p\in\Nset$ (thus $t_p=pT$ here).
As in Example \ref{ex2}, the inequalities here are meant componentwise.

The information transfers are schematised on Figure \ref{fi4}.
The agents are numbered with Arabic numbers and the Roman numbers describe the different stages of the spanning completion.
Only the communications with guaranteed coefficient $\gamma$ are represented.
For simplicity, the self-loops are omitted.
\begin{figure}[h]
\begin{center}
\input fig4.pstex_t
\caption{see text of Example \ref{ex55}.}
\label{fi4}
\end{center}
\end{figure}

Analysing the system with the use of a unique spanning-tree (Theorem \ref{th1}) yields: $\trho_1= \left(
1-\gamma^3
\right)^{1/5}$.

To use Theorem \ref{th3} for analysis, one considers two spanning-trees and takes
on any interval $[t_p,t_{p+1}]$:
\begin{gather*}
\cN^1_{p,0}=\cN^2_{p,0}=\{1\},\quad
\cN^1_{p,1}=\cN^2_{p,1}=\{1,2\},\\
\cN^1_{p,2}=\{1,2,3\},\
\cN^2_{p,2}=\{1,2,4\},\quad
\cN^1_{p,3}=\{1,2,3,5\},\
\cN^2_{p,3}=\{1,2,4,6\},\\
\cN^1_{p,4}=\{1,2,3,5,6\},\
\cN^2_{p,4}=\{1,2,4,5,6\},\quad
\cN^1_{p,5}=\cN^2_{p,5}=\cN
\end{gather*}
and
\[
A(0)=\dots=A(4)=\begin{pmatrix}
\gamma & 0\\ 0 & \gamma
\end{pmatrix}\ .
\]
The deduced estimate is $\trho_2= \left(1-\gamma^5\right)^{1/5}$.
The important point is that it is systematically {\em looser} than the previous one.
Indeed, the previous formula could have been obtained by taking into account only one of the two spanning-trees, say the ``right branch", where the signal circulates in the order 1--2--4--6--5--3.
Otherwise said, there would no difference in evaluating the graph similarly schematised, shown on Figure \ref{fi5}.
\begin{figure}[h]
\begin{center}
\input fig5.pstex_t
\caption{see text of Example \ref{ex55}.}
\label{fi5}
\end{center}
\end{figure}

How to take into account the crossing of the two spanning-trees, a case explicitly discarded in Section \ref{se6}?
A general idea is to introduce new ``populations".
However, this is not so easy, as Lemma \ref{le1} is hardly adapted
(it is here useful to recall that the diameter of the union of two sets is at most equal to the sum of the diameters of these two sets, {\em if their intersection is non-void}).

Along these lines, one may propose an idea for improvement of $\trho_2$.
Let
\[
\cN^3_{p,t}\doteq\cN^1_{p,t}\cup\cN^2_{p,t}\ .
\]
This is just, in fact, the population considered in the one-spanning-tree method leading to $\trho_1$.
We are then allowed to take:
\begin{equation}
\label{A55}
A(0)=\dots=A(2)=
\begin{pmatrix}
\gamma & 0 & 0\\ 0 & \gamma & 0\\ 0 & 0 & \gamma
\end{pmatrix},\
A(3)=A(4)=
\begin{pmatrix}
\gamma & 0 & 1-\gamma\\ 0 & \gamma & 1-\gamma\\ 0 & 0 & \gamma
\end{pmatrix}\ .
\end{equation}
As
\[
A(4)A(3)A(2)A(1)A(0)=\gamma^4\begin{pmatrix}
\gamma & 0 & 2(1-\gamma)\\
0 & \gamma & 2(1-\gamma)\\
0 & 0 & \gamma
\end{pmatrix}\ ,
\]
the estimate obtained via Theorem \ref{th3} is
\[
\trho_3\doteq\left(
1-\gamma^4(2-\gamma)
\right)^{1/5}\ ,
\]
which verifies $\trho_1\leq\trho_3\leq\trho_2$ for $\gamma\in [0,1]$: actually, $\trho_3$  does not overpass the precision of $\trho_1$.
\qed
\end{exam}

A careful examination of the previous example shows why no improvement could be obtained: the diameters of the three sets are equal up to the third stage, and the form of the difference inequalities involved forbid the two components fed with by the third one, to become larger than the latter.

However, notice that this paradoxical behaviour is also resulting of the value of the coefficients.
The next example indicates that the method proposed in Example \ref{ex55} can indeed provide better estimates.

\begin{exam}
\label{ex8}
We consider a slight modification of Example \ref{ex55}.
For fixed $\eta\in [0,1]$, we take $\Gamma$ as previously (see \eqref{G55}), except
\[
\Gamma(pT+1)
\geq\begin{pmatrix}
\gamma & 0 & 0 & 0 & 0 & 0\\
0 & \gamma & 0 & 0 & 0 & 0\\
0 & \eta\gamma & 0 & 0 & 0 & 0\\
0 & \gamma & 0 & 0 & 0 & 0\\
0 & 0 & 0 & 0 & 0 & 0\\
0 & 0 & 0 & 0 & 0 & 0
\end{pmatrix},\
\Gamma(pT+2)
\geq\begin{pmatrix}
\gamma & 0 & 0 & 0 & 0 & 0\\
0 & \gamma & 0 & 0 & 0 & 0\\
0 & 0 & \eta\gamma & 0 & 0 & 0\\
0 & 0 & 0 & \gamma & 0 & 0\\
0 & 0 & \gamma & 0 & 0 & 0\\
0 & 0 & 0 & \gamma & 0 & 0
\end{pmatrix}\ .
\]
In other words, the transmission along the ``left branch" in Figure \ref{fi4} occurs with a least coefficient than along the right one.
This modifies both the evolution of the diameters of $\cN^1$ and $\cN^3$, and one now has to modify the values of $A$ by taking
\[
A(1)=A(2)=
\begin{pmatrix}
\eta\gamma & 0 & 0\\ 0 & \gamma & 0\\ 0 & 0 & \eta\gamma
\end{pmatrix}\ ,
\]
instead of those given in \eqref{A55}.
Using the notations of Example \ref{ex55} yields the two contraction rate estimates
\[
\trho_1'\doteq\left(
1-\eta^2\gamma^3
\right)^{1/5}\quad
\text{ and }\quad
\trho_3'\doteq\left(
1-\gamma^4(\gamma+2\eta^2(1-\gamma)
\right)^{1/5}\ .
\]
In particular, when
\[
\eta\leq
\frac{\gamma}{(1+2\gamma(1-\gamma))^{1/2}}
\]
(a quantity located in $[0,1/3]$ for $\gamma\in [0,1/2]$), then the $\trho_3'$ is smaller than the estimate $\trho_1'$, obtained by considering a single spanning-tree.
\qed
\end{exam}
 
\section{Conclusion}
\label{se77}
Several tools for estimating the convergence rate to consensus in multiagent systems were introduced and
illustrated through simple examples. The criteria are based on topological as well as basic quantitative
information.
In accordance to previous results, consensus is reached provided that information can flow at least along some
spanning tree from one agent to all of the others.
A key quantity, in this respect, appears to be a lower bound on the total weight of the agents located upstream along the information flow for any chosen spanning-tree.
More general criteria are also provided in which tighter estimates are allowed, provided that more spanning-trees
are simultaneously taken into account.

These techniques are, in general, based on the idea of considering a decomposition of the overall population into
subsets which influence each other in some quantifiable ways. Natural candidates for this partition
appear to be the agents already attained by the information flows along the spanning-trees. 
There seem to be technical difficulties in trying to consider other kinds of partitions, as in general, neither the diameter of an union of sets, nor of an intersection of sets, is related to the diameters of the two sets.
However, it may be possible to consider the set of agents attained by one ore more spanning-trees and then  the set of agents attained in the reverse order. We leave this as an interesting open question for future research.

The method presented here provides results which are rather tight and inherently robust due to the qualitative nature of the assumptions involved.
It is especially interesting to develop tools for quantitative estimates based on the consideration of simultaneous trees as arising from a single tree which gets repeated through time, as in  Example \ref{ex6}.
Again this will be topic of further investigations.

\appendix

\section{Appendix -- Fundamental inequalities}


We state in the sequel a result on difference inequalities which is central to the techniques developed in the text.
Consider the time-varying linear system \eqref{equb}.
As before, the index set $\cN$ is finite or countable, the $x_k$ constitute a collection of scalar functions defined on $\Nset$, $\cI$ is a finite or countable index set and, for any $t\in\Nset$, a collection of subsets $\cM_i(t)$ of $\cN$, $i\in\cI$, is given.
Also, the state-matrices $(\gamma_{k,l}(t))_{(k,l)\in\cN\times\cN}$ of the system are row-stochastic.
Define the diameters:
\begin{equation*}
\Delta(t)\doteq \sup_{k\in\cN} x_k(t) - \inf_{k\in\cN}x_k(t),\qquad
\Delta_i(t)\doteq \diam\ \cM_i(t)=\sup_{k\in\cM_i(t)} x_k(t) - \inf_{k\in\cM_i(t)}x_k(t),
\end{equation*}
and the vector
\begin{equation*}
\Delta_\cM(t) \doteq (\diam\ \cM_i(t))_{i\in\cI} = (\Delta_i(t))_{i\in\cI}\ .
\end{equation*}
The following result provides informations on the evolution of the diameter vector.

\begin{lemm}
\label{le1}
Assume that for all $k\in\cN$, for all $t\in\Nset$, some sets $\cM_{k,j}(t)$, $j\in\cI$, are given, such that
\[
\cM_{k,j}(t)\subseteq\cM_j(t)\quad\text{ and }\quad
\cM_{k,j}(t)\cap\cM_{k,j'}(t)\neq\emptyset\Rightarrow j=j'\ .
\]
Let maps $c_{i,j}(t)$, $i,j\in\cI$, and $C(t)$ be such that:
\begin{equation}
\label{aA}
c_{i,j}(t) \leq \inf_{k\in\cM_i(t+1)} \sum_{l\in\cM_{k,j}(t)} \gamma_{k,l}(t),\qquad
C(t)\doteq (c_{i,j}(t))_{(i,j)\in\cI\times\cI}\ .
\end{equation}

Then, for any $t,T\in\Nset$, for any $i\in\cI$,
\begin{equation}
\label{ineq}
0\leq
\Delta_i(t+T)
\leq
\bfo_i\t \left(
\vphantom{\int}
C(t+T-1)\dots C(t) \Delta_\cM(t)
+\left(
1-C(t+T-1)\dots C(t)
\right)
\bfo\Delta(t)\right)\ .
\end{equation}
\end{lemm}

By convention, we put
$\sum_{i \in \emptyset} c_i =0$ and $\inf_{i \in \emptyset} c_i = + \infty$.
Recall that the vector $\bfo$ in the statement is made up of a column of $1$ and that the vector $\bfo_i$ has null components, except 1 in the $i$-th position (in finite dimension, it is the $i$-th vector of the canonic basis).
In particular, $\bfo_i\t \Delta_\cM(t)=\Delta_i(t)$, $\bfo=\sum_{i\in\cI}\bfo_i$.

\begin{rema}
Notice that formula \eqref{ineq} may involve infinite summations in the products of infinite-dimen\-sio\-nal matrices.
As the coefficients of the matrices $C(t)$ are nonnegative and bounded by 1, uniform convergence of the series of terms indeed occurs on any bounded time interval, therefore the notation has a univocal meaning.
\qed
\end{rema}

\begin{proof}
Define
\[
M(t)\doteq \sup_{k\in\cN} x_k(t),\
M_i(t)\doteq \sup_{k\in\cM_i(t)} x_k(t),\
m(t)\doteq \min_{k\in\cN}x_k(t),\
m_i(t)\doteq \min_{k\in\cM_i(t)}x_k(t)\ ,
\]
in such a way that the quantities previously defined in the statement verify:
\[
\Delta\equiv M-m,\
\Delta_i\equiv M_i-m_i\ .
\]

First of all, notice that, due to the nonnegativity of the coefficients $\gamma_{k,l}(t)$, identity \eqref{equb} implies, for any $t\in\Nset$ and for any $k\in\cN$,
\[
x_k(t+1) \leq \sum_{l\in\cN} \gamma_{k,l}(t) M(t) =  M(t)\ .
\]
Taking the supremum and arguing similarly for the lower bounds, we obtain:
\begin{equation*}
M(t+1)\leq M(t),\ m(t+1) \geq m(t)\ .
\end{equation*}
In particular,
\begin{equation}
\label{cN1b}
\Delta(t+1)\leq \Delta(t)\ .
\end{equation}

Also, due to the fact that $\cM_i(t)\subseteq\cN$, it comes:
\begin{equation}
\label{cN2a}
M_i(t)= \sup_{k\in\cM_i(t)} x_k(t)\leq \sup_{k\in\cN} x_k(t) = M(t),\quad
m_i(t) \geq m(t)\ ,
\end{equation}
and
\begin{equation*}
\Delta_i(t) \leq \Delta(t)\ .
\end{equation*}

Applying tighter estimate, one obtains from \eqref{equb} that, for any $k\in\cN$,
\begin{eqnarray*}
x_k(t+1)
& = &
\sum_{j\in\cI}\sum_{l\in\cM_{k,j}(t)}\gamma_{k,l}(t)x_l(t)
+\sum_{l\in\cN\setminus\bigcup_{j\in\cI}\cM_{k,j}(t)}\gamma_{k,l}(t) x_l(t)\\
& \leq &
\sum_{j\in\cI} \left(
\sum_{l\in\cM_{k,j}(t)}\gamma_{k,l}(t)
\right) M_j(t)
+\sum_{l\in\cN\setminus\bigcup_{j\in\cI}\cM_{k,j}(t)}\gamma_{k,l}(t) M(t)\\
& = &
\sum_{j\in\cI}\sum_{l\in\cM_{k,j}(t)}\gamma_{k,l}(t)M_j(t)
+\left(
1-\sum_{j\in\cI}\sum_{l\in\cM_{k,j}(t)}\gamma_{k,l}(t)
\right)M(t)\ ,
\end{eqnarray*}
due to \eqref{aA} and the nonnegativity of the coefficients $\gamma_{k,l}(t)$.
 If now $k\in\cM_i(t+1)$ for some $i\in\cI$, one obtains:
\begin{eqnarray*}
x_k(t+1)
& \leq &
\sum_{j\in\cI}c_{i,j}(t)M_j(t)
+\left(
1-\sum_{j\in\cI}c_{i,j}(t)
\right)M(t)
+\sum_{j\in\cI}\left(
c_{i,j}(t)-\sum_{l\in\cM_{k,j}(t)}\gamma_{k,l}(t)
\right)
\left(
M(t)-M_j(t)
\right)\\
& \leq &
\sum_{j\in\cI}c_{i,j}(t)M_j(t)
+\left(
1-\sum_{j\in\cI}c_{i,j}(t)
\right)M(t)\ ,
\end{eqnarray*}
due to the fact that $M(t)\geq M_j(t)$ for any $t\in\Nset$ and any $j\in\cI$ (see \eqref{cN2a}).
Consequently, for any $i\in\cI$,
\[
M_i(t+1)\leq \sum_{j\in\cI}c_{i,j}(t)M_j(t)
+\left(
1-\sum_{j\in\cI}c_{i,j}(t)
\right)M(t)\ .
\]

One establishes similarly that
\[
x_k(t+1)
\geq
\sum_{j\in\cI}c_{i,j}(t)m_j(t)
+\left(
1-\sum_{j\in\cI}c_{i,j}(t)
\right)m(t)\ ,
\]
with the {\em same} coefficients, so, for any $i\in\cI$:
\[
m_i(t+1)
\geq
\sum_{j\in\cI}c_{i,j}(t)m_j(t)
+\left(
1-\sum_{j\in\cI}c_{i,j}(t)
\right)m(t)\ .
\]

Subtracting the previous inequalities, one may thus deduce that, for any $i\in\cI$,
\begin{equation*}
\Delta_i(t+1)
\leq
\sum_{j\in\cI}c_{i,j}(t)\Delta_j(t)
+\left(
1-\sum_{j\in\cI}c_{i,j}(t)
\right)\Delta(t)\ .
\end{equation*}

The collection of these inequalities, together with \eqref{cN1b}, may be written under the matrix (possibly infinite) form:
\[
\begin{pmatrix}
\Delta_\cM(t+1)\\ \Delta(t+1)
\end{pmatrix}
\leq
\begin{pmatrix}
C(t) & \bfo -C(t)\bfo\\
0 & 1
\end{pmatrix}
\begin{pmatrix}
\Delta_\cM(t)\\ \Delta(t)
\end{pmatrix}\ .
\]
The previous inequality has to be understood componentwise.

Now, one shows easily that:
\[
\begin{pmatrix}
C(t+1) & \bfo -C(t+1)\bfo\\
0 & 1
\end{pmatrix}
\begin{pmatrix}
C(t) & \bfo -C(t)\bfo\\
0 & 1
\end{pmatrix}
= \begin{pmatrix}
C(t+1)C(t) & \bfo -C(t+1)C(t)\bfo\\
0 & 1
\end{pmatrix}\ ,
\]in such a way that, for nonnegative $T$,
\[
\begin{pmatrix}
\Delta_\cM(t+T)\\ \Delta(t+T)
\end{pmatrix}
\leq
\begin{pmatrix}
C(t+T-1)\dots C(t) & \bfo -C(t+T-1)\dots C(t)\bfo\\
0 & 1
\end{pmatrix}
\begin{pmatrix}
\Delta_\cM(t)\\ \Delta(t)
\end{pmatrix}\ .
\]
This formula permits to complete the proof of Lemma \ref{le1}.
\end{proof}

\section{Appendix -- Proof of Lemma \protect{\ref{le4}}}

\noindent {\bf 1.}
One verifies directly that, for any $t=0,\dots, n+q-3$,
\[
A(t) =  I_q-(1-\gamma) \left(
\vphantom{\int}
\delta_{0\leq t\leq n-1}\ e_1e_1\t +\delta_{1\leq t\leq n}\ e_2(e_2-e_1)\t +\dots+\delta_{q-1\leq t\leq n+q-2}\ e_q(e_q-e_{q-1})\t
\right)\ ,
\]
where $e_i$ is the $i$-th vector of the canonical basis in $\Rset^q$ and $\delta_{i\leq t \leq i+n-1}$ is 1 (resp.\ 0) if the condition written in index is fulfilled (resp.\ violated).

Let us first establish the following factorisation formula:
\begin{equation}
\label{gene}
A(n+q-3)\dots A(0)
=
\gamma^{n-q}\ \diag\{ 1;\gamma;\dots; \gamma^{q-1}\}\
B(n+q-3)\dots B(0)\
\diag\{ \gamma^{q-1};\dots ;\gamma;1\}\ ,
\end{equation}
where the matrix $B(t)$ is obtained from $A(t)$ by replacing $\gamma$ on the diagonal by 1, and $1-\gamma$ by
\[
\xi\doteq\frac{1-\gamma}{\gamma}\ ,
\]
that is simply:
\[
B(t)\doteq
I_q+\xi \sum_{i=1}^{q-1}\delta_{i\leq t\leq i+n-2}\ e_{i+1}e_i\t\ .
\]
In formula \eqref{gene} and below, $\diag$ is used to define diagonal matrices.

Formula \eqref{gene} will be proved by induction on the positive integer $q$.
Notice that strictly speaking, the matrices $A, B\in\Rset^{q\times q}$ depend upon $q$ (and $n$), but for simplicity we omit here any explicit indication of this dependence.
Indeed, for $q=1$, $A(t)=\gamma$ for $0\leq t\leq n+q-3=n-2$, and $A(n-2)\dots A(0)=\gamma^{n-1}$; while for $q=2$, $n+q-3=n-1$ and
\[
A(0) = \begin{pmatrix}
\gamma & 0\\ 0 & 1
\end{pmatrix},\quad
A(t)=\begin{pmatrix}
\gamma & 0\\ 1-\gamma & \gamma
\end{pmatrix},\ 1\leq t\leq n-2,\quad
A(n-1)=\begin{pmatrix}
1 & 0\\ 1-\gamma & \gamma
\end{pmatrix}\ ,
\]
so that
\[
A(n-1)\dots A(0)=\gamma^{n-2}
\begin{pmatrix}
1 & 0\\ 0 & \gamma
\end{pmatrix}
\begin{pmatrix}
1 & 0\\ \xi & 1
\end{pmatrix}^{n-1}
\begin{pmatrix}
\gamma & 0\\ 0 & 1
\end{pmatrix}
=\gamma^{n-2}
\begin{pmatrix}
1 & 0\\ 0 & \gamma
\end{pmatrix}
B(n-1)\dots B(0)
\begin{pmatrix}
\gamma & 0\\ 0 & 1
\end{pmatrix} \ .
\]
(Notice that $B(0)=I_2$ and $B(t)=\begin{pmatrix}
1 & 0\\ \xi & 1
\end{pmatrix}$ for $t=1,\dots, n-1$).

Assume now that \eqref{gene} is true at order $q-1$ and consider order $q$.
Due to the particular structure of the matrices $A$ and $B$, which are null except terms on the diagonal and the sub-diagonal, one has
\begin{subequations}
\label{comm}
\begin{gather}
\label{comma}
\diag\{I_{q-1}; 0\}
\left(
\prod_{t=0}^{n+q-3} A(t)
\right)
\diag\{I_{q-1}; 0\}
=
\prod_{t=0}^{n+q-3}
\diag\{I_{q-1}; 0\}
A(t)\
\diag\{I_q; 0\},\\
\label{commb}
\diag\{0; I_{q-1}\}
\left(
\prod_{t=0}^{n+q-3} A(t)
\right)
\diag\{0; I_{q-1}\}
=
\prod_{t=0}^{n+q-3}
\diag\{0; I_{q-1}\}
A(t)\
\diag\{0; I_{q-1}\}\ ,
\end{gather}
\end{subequations}
and similarly for $B(t)$.
In the previous identities and in subsequent formulas, the products are non-commutative: the convention is that $t$ is decreasing from the left factor to the right one.

Now, it is easy to identify the right-hand sides of the two identities \eqref{comm} with a product of $n+(q-1)-3=n+q-4$ matrices $A$ (resp.\ $B$) corresponding to the index $q-1$ (the last term in the right-hand product in \eqref{comma}, resp.\ the first term in the right-hand product in \eqref{commb}, is equal to $\diag\{I_{q-1}; 0\}$, resp.\ $\diag\{0; I_{q-1}\}$, and can be suppressed).
Using the induction hypothesis at order $q-1$, one shows that
\begin{eqnarray*}
\lefteqn{\diag\{I_{q-1}; 0\}
\left(
\prod_{t=0}^{n+q-3} A(t)
\right)
\diag\{I_{q-1}; 0\}}\\
& = &
\gamma^{n-(q-1)}\ 
\diag\{ 1;\dots ;\gamma^{q-2}; 0\}\
\diag\{I_{q-1}; 0\}\
\left(
\prod_{t=0}^{n+q-3} B(t)
\right)\
\diag\{I_{q-1}; 0\}\
\diag\{ \gamma^{q-2};\dots ;1;0\}\\
& = &
\gamma^{n-(q-1)}\gamma^{-1}\
\diag\{I_{q-1}; 0\}\
\diag\{ 1;\dots ; \gamma^{q-1}\}\
\left(
\prod_{t=0}^{n+q-3} B(t)
\right)\
\diag\{ \gamma^{q-1};\dots ;1\}\
\diag\{I_{q-1}; 0\}\\
& = &
\gamma^{n-q}\
\diag\{I_{q-1}; 0\}\
\diag\{ 1;\dots; \gamma^{q-1}\}\
\left(
\prod_{t=0}^{n+q-3} B(t)
\right)\
\diag\{ \gamma^{q-1};\dots ;1\}\
\diag\{I_{q-1}; 0\}\ .
\end{eqnarray*}
One establishes similarly that
\begin{multline*}
\diag\{0;I_{q-1}\}
\left(
\prod_{t=0}^{n+q-3} A(t)
\right)
\diag\{0;I_{q-1}\}\\
=
\gamma^{n-q}\
\diag\{0;I_{q-1}\}\
\diag\{ 1;\dots; \gamma^{q-1}\}\
\left(
\prod_{t=0}^{n+q-3} B(t)
\right)\
\diag\{ \gamma^{q-1};\dots ;1\}\
\diag\{0;I_{q-1}\}\ ,
\end{multline*}
and this is indeed sufficient, due to the structure of the matrices $A$ and $B$ mentioned earlier, to prove that \eqref{gene} is true at order $q$.
This achieves the proof of \eqref{gene} by induction.\\

\noindent {\bf 2.}
One now estimates the matrix-product
\[
\Pi=(\Pi_{i,j})_{(i,j)\in\{1,\dots,q\}^2}\doteq
\prod_{t=0}^{n+q-3} B(t)
=
\prod_{t=1}^{n+q-3}
\left(
I_q+\xi \sum_{i=1}^{q-1}\delta_{i\leq t\leq i+n-2}\ e_{i+1}e_i\t
\right)\ .
\]
Each term of this product is a lower-triangular matrix, so $\Pi$ shares the same property.

The fact that the canonical basis is orthonormal implies that, for any $i> j$, $i,j\in\{1,\dots, q\}$, it holds:
\[
\Pi_{i,j}\doteq \xi^{i-j}\ \card\left\{
\vphantom{\int}
(t_{j+1},\dots, t_i)\in [j,j+n-2]\times\dots\times [i-1,i+n-3]\cap\Nset^{i-j}\ :\  t_{j+1}<\dots <t_i 
\right\}\ ,
\]
and also that the diagonal terms are equal to 1.
The previous formula just means that, for a term in $e_ie_j\t$ to emerge from the product, it should be the result of the product
\[
(e_ie_{i-1}\t) \cdot (e_{i-1}e_{i-2}\t) \dots (e_{j+1}e_j\t)\ ,
\]
where each of the term between parentheses comes from a certain matrix $A(t)$ --- the rest of the factors coming from identity matrices.
Conversely, all products of different type vanishes.

In order to evaluate the quantities $\Pi_{i,j}$ previously defined, notice that the change of variables
\[
t'_{j+1}=t_{j+1},\ t'_{j+2}=t_{j+2}-1,\ \dots,\ t'_i=t_i-(i-j-1)
\]
yields:
\begin{multline*}
\card\left\{
\vphantom{\int}
(t_{j+1},\dots, t_i)\in [j,j+n-2]\times\dots\times [i-1,i+n-3]\cap\Nset^{i-j}\ :\  t_{j+1}<\dots <t_i 
\right\}\\
=\card\left\{
\vphantom{\int}
(t'_{j+1},\dots, t'_i)\in \left(
[j,j+n-2]\cap\Nset
\right)^{i-j}\ :\  t'_{j+1}\leq\dots\leq t'_i 
\right\}\ .
\end{multline*}

\noindent {\bf 3.}
We now compute explicitly the value of the function $F(m,n)$ defined on $\Nset\times\Nset$ as:
\[
F(m,n)\doteq\card\left\{
\vphantom{\int}
(t_1,\dots, t_m)\in \left(
[1,n]\cap\Nset
\right)^m\ :\  t_1\leq\dots\leq t_m 
\right\}\ .
\]
Clearly,
\[
F(1,n)=n,\quad F(2,n)=\frac{n(n+1)}{2}\ .
\]

Considering separately the cases where $t_1=1$, $t_1=2$, \dots, $t_1=n$, one finds the following induction relation:
\[
F(m,n)=\sum_{i=1}^n F(m-1,i)\ .
\]
On the other hand, let
\[
G(m,n)\doteq\begin{pmatrix}
m+n-1\\ m
\end{pmatrix}
=\frac{(m+n-1)!}{m! (n-1)!} \ ,
\]
one has:
\[
G(1,n)=\begin{pmatrix}
n\\ 1
\end{pmatrix} =n,\quad G(2,n)=\begin{pmatrix}
n+1\\ 2
\end{pmatrix}=\frac{n(n+1)}{2}\ .
\]
Independently, it is known that
\[
\begin{pmatrix}
n\\ m
\end{pmatrix}
=
\begin{pmatrix}
n-1\\ m-1
\end{pmatrix}
+ \begin{pmatrix}
n-1\\ m
\end{pmatrix}\ ,
\]
in such a way that
\[
G(m,n)=\begin{pmatrix}
m+n-1\\ m
\end{pmatrix}=
\begin{pmatrix}
m+n-2\\ m-1
\end{pmatrix}
+\begin{pmatrix}
m+n-2\\ m
\end{pmatrix}
=G(m-1,n)+G(m,n-1)\ .
\]
It ensues from repeated use of this formula, that:
\begin{multline*}
G(m,n)=G(m-1,n)+G(m,n-1)
=G(m-1,n)+G(m-1,n-2)+G(m,n-3)\\
=\dots
= \sum_{i=2}^n G(m-1,i)+G(m,1)
=\sum_{i=1}^n G(m-1,i)\ ,
\end{multline*}
because $G(m,1)=G(m-1,1)=1$.
Having the same initial condition and sharing the same induction relation, $F$ and $G$ are thus equal, and $F(m,n)=\begin{pmatrix}
m+n-1\\ m
\end{pmatrix}$.\\

\noindent {\bf 4.}
The value of $F$ found before is now used to estimate $\Pi$ and next $\trho_q$.
We deduce from what precedes that, for $i>j$,
\[
\Pi_{i,j}=\xi^{i-j} F(i-j,n-1)=\xi^{i-j}\begin{pmatrix}
i-j+n-2 \\ i-j
\end{pmatrix}
=\xi^{i-j}\begin{pmatrix}
i-j+n-2 \\ n-2
\end{pmatrix}\ .
\]
Recall that $\Pi_{i,i}=1$ and $\Pi_{i,j}=0$ for $i<j$.

From the fact that the matrix $\Pi$ above is lower-triangular, one finds out by application of Theorem \ref{th3} that 
\[
1-\trho_q^{n+q-2}
= \max_{i=1,\dots, q} \bfo_iA(n+q-3)\dots A(0)\bfo
=\bfo_qA(n+q-3)\dots A(0)\bfo\ .
\]
From \eqref{gene} and the previous computations, one thus deduces 
\begin{multline*}
1-\trho_q^{n+q-2} = \gamma^{n-1}\sum_{i=1}^q
\gamma^{q-i}\Pi_{q,i}
=\gamma^{n-1}\sum_{i=1}^q
\gamma^{q-i}\xi^{q-i}\begin{pmatrix}
q-i+n-2 \\ n-2
\end{pmatrix}\\
= \gamma^{n-1}\sum_{i=1}^q
\gamma^{q-i}\left(
\frac{1-\gamma}{\gamma}
\right)^{q-i}\begin{pmatrix}
q-i+n-2 \\ n-2
\end{pmatrix}\ .
\end{multline*}
Thus,
\[
\trho_q^{n+q-2}
= 1-\gamma^{n-1}
\sum_{i=1}^q
(1-\gamma)^{q-i}\begin{pmatrix}
q-i+n-2 \\ n-2
\end{pmatrix}\ ,
\]
or again:
\[
\trho_q
= \left(
1- \gamma^{n-1}
\sum_{i=0}^{q-1}
(1-\gamma)^i\begin{pmatrix}
n+i-2 \\ n-2
\end{pmatrix}
\right)^{1/(n+q-2)}\ .
\]
This achieves the proof of the first equality in the statement of Lemma \ref{le4}.\\

\noindent {\bf 5.}
To show the identity of the two expressions in Lemma \ref{le4}, notice that
\begin{multline*}
\sum_{i=0}^{q-1} \begin{pmatrix}
n+i-2\\ n-2
\end{pmatrix}\delta^i
= \frac{1}{(n-2)!}\sum_{i=0}^{q-1} 
\frac{d^{n-2}}{d\delta^{n-2}} \left[
\delta^{n+i-2}
\right]
= \frac{1}{(n-2)!}
\frac{d^{n-2}}{d\delta^{n-2}} \left[
\sum_{i=0}^{q-1} \delta^{n+i-2}
\right]\\
= \frac{1}{(n-2)!}
\frac{d^{n-2}}{d\delta^{n-2}} \left[
\frac{1-\delta^{n+q-2}}{1-\delta}
\right]\ .
\end{multline*}
On the other hand, one shows easily that
\[
\frac{1}{(n-2)!}
\frac{d^{n-2}}{d\delta^{n-2}} \left[
\frac{1}{1-\delta}
\right]
= \frac{1}{(1-\delta)^{n-1}}\ .
\]
This permits to deduce the identity of the two expressions in the statement.\\
 
 \noindent {\bf 6.}
Last, we show the limiting property expressed in Lemma \ref{le4}.
From the last formula, one may see that, for every $n\geq 2$:
\[
\trho_q = \sqrt[n+q-2]{ P_n (q, \delta) \delta^{q} }\ ,
\]
where $P_n$ is a polynomial in $q$ and $\delta=1-\gamma$ of degree $n-2$ with respect to both variables.
Henceforth, taking the limit for $q \rightarrow + \infty$ yields the estimate:
\[
\lim_{q \rightarrow + \infty}\trho_q
= \lim_{q \rightarrow + \infty} \sqrt[n+q-2]{ P_n (q, \delta) \delta^{q} } = \delta =1-\gamma\ ,
\]
which corresponds to the true value of the converging rate.
Indeed,
\[
\sqrt[n+q-2]{ P_n (q, \delta) }
= e^{\ln P_n (q, \delta)/(n+q-2)}
= e^{\left[(n-2)\ln q +\ln (1+O(1/q))\right]/(n+q-2)}\ ,
\]
as $P_n$ is of degree $n-2$ in $q$.
The asymptotic expansion announced in the statement is thus proved, and this achieves the proof of Lemma \ref{le4}.

\end{document}